\newcommand{\MailFile}[1]{}
\def\PaperDate{2008/08/14}
\documentclass[12pt,letterpaper]{article}
  \usepackage{buxmath}
  \usepackage{bux2ref}
  \usepackage[german,american]{babel}
  \usepackage{color}
  \usepackage{ulem}
  \usepackage{datum}

  \usepackage{email}
  \MailFile{rank_one_files.tex}
  \DoNotMail{ulasy.fd,babel.*,*.ldf,ams*.sty,*.bref}

  \newcommand{\squish}[2][-2pt]{%
    \begingroup
      \setlength{\fboxsep}{#1}%
      \setlength{\fboxrule}{0pt}%
      \kern-\fboxsep
      \framebox{\ensuremath{#2}}%
      \kern-\fboxsep
    \endgroup
  }

\iftrue
  \DeclareRobustCommand{\mydelete}[1]{%
    \begin{color}{green}\sout{#1}\end{color}
  }
  \DeclareRobustCommand{\myinsert}[1]{\begin{color}{red}#1\end{color}}
  \DeclareRobustCommand{\myswap}[2]{\sout{#1} \myinsert{#2}}
  \DeclareRobustCommand{\change}[2][]{%
    \begingroup
      \def\next{}%
      \def\deleted{#1}%
      \def\newmat{#2}%
      \def\empty{}%
      \ifx\deleted\empty
        \def\next{\myinsert{#2}}%
      \else
        \ifx\newmat\empty
          \def\next{\mydelete{#1}}%
        \else
          \def\next{\myswap{#1}{#2}}%
        \fi
      \fi
    \expandafter\endgroup\next
  }
\else
  \DeclareRobustCommand{\change}[2][]{%
    \begingroup
      \def\next{}%
      \def\newmat{#2}%
      \def\empty{}%
      \ifx\newmat\empty
        \def\next{\unskip}%
      \else
        \def\next{#2}%
      \fi
    \expandafter\endgroup\next
  }
\fi

  \InputIfFileExists{basic_notation.tex}{}{}

  \newvariable{\TheIdentityElement}{1}
  \newvariable{\ThePrime}{v}
  \newvariable{\ThePrimeSet}{S}
  \newvariable{\OkaRing}{\Oka_S}
  \newvariable{\TheGroup}{\mathbf{G}}
  \newvariable{\AltGroup}{\mathbf{H}}
  \newvariable{\TheTorus}{\mathbf{A}}
  \newvariable{\TheTorusElement}{a}
  \newvariable{\TheParabolic}{\mathbf{P}}
  \newvariable{\TheRadical}{\mathbf{R}^{\mathbf{u}}}

  \newvariable{\TheUnipotent}{\mathbf{U}}
  \newvariable{\TheLevi}{\mathbf{L}}
  \newvariable{\AltTorus}{\mathbf{A}^0}

  \newvariable{\TheConjugateUnip}{u}
  \newvariable{\TheCharacter}{\chi}
  \newvariable{\TheField}{K}
  \newvariable{\SumOfLocalRanksFct}{k}
  \newvariable{\SumOfLocalRanks}{m}

  \newvariable{\TheCutOff}{q}
  \newvariable{\TheCompactSet}{C}
  \newvariable{\TheFiniteSet}{F}
  \newvariable{\TheFiniteSetElement}{f}
  \newvariable{\AltFiniteSetElement}{\hat{f}}

  \newvariable{\TheBaseVertex}{e}
  \newvariable{\TheRay}{\rho}
  \newvariable{\TheParameter}{t}
  \newvariable{\ThePoint}{x}
  \newvariable{\AltPoint}{y}
  \newvariable{\TheHoroRadius}{r}
  \newvariable{\TheThickness}{L}
  \newvariable{\TheCoverDepth}{T}
  \newvariable{\TheSliceSize}{R}
  \newvariable{\TheHoroDepth}{R}
  \newvariable{\ThePush}{d}

  \newvariable{\TheLattice}{\Gamma}
  \newvariable{\TheLatticeElement}{\gamma}
  \newvariable{\TheArithmeticElement}{\gamma}
  \newvariable{\TheRadius}{r}
  \newvariable{\TheDistance}{M}
  \newvariable{\TheLag}{L}
  \newvariable{\TheExponent}{n}
  \newvariable{\TheBigExponent}{N}

  \newcommand{\ParenthesesOf}[2]{\left(#1\right)\left(#2\right)}

  \newcommand{\DummyArg}{\,\,\cdot\,\,}
  \newcommand{\NormOf}[2][]{\left|#2\right|_{#1}}
  \newcommand{\acts}{\cdot}

  \newvariable{\TheGroupElement}{g}

  \newvariable{\TheHoroCollection}{\mathcal{H}}

  \newvariable{\FiniteField}{\FFF}
  \newvariable{\ThePrimePower}{q}
  \newvariable{\TheIndeterminate}{t}

  \newvariable{\AbstractGroup}{\Gamma}

  \newvariable{\Sl}{\operatorname{SL}}

  \newcommand{\pref}[1]{(\ref{#1})}

  \newcommand{\isom}{\cong}
  \newcommand{\subgroup}{\leq}

  \newcommand{\AbsValueOf}[1]{|#1|}

  \newcommand{\InteriorOf}[1]{\topcirc{#1}}

  \newcommand{\RankOf}[2][]{\operatorname{rank}_{#1}(#2)}

  \newvariable{\TheType}{m}
  \newcommand{\FType}[1]{F${}_{#1}$}
  \newcommand{\EilenbergMaclaneOf}[1]{K(#1,1)}

  \newvariable{\AltPecComplex}{Y}

  \newvariable{\TheConnectivity}{m}
  \newvariable{\AltConnectivity}{n}

  \newvariable{\TheReal}{s}
  \newvariable{\AltReal}{t}
  \newcommand{\Pair}[2]{(#1,#2)}
  \newcommand{\crossprod}{\times}

  \newvariable{\Link}{\operatorname{Lk}}
  \newvariable{\DescLink}{\operatorname{Lk}^{\downarrow}}
  \newvariable{\EquLink}{\operatorname{Lk}^{=}}
  \newvariable{\VertLink}{\operatorname{Lk}_{\ver}}
  \newvariable{\HorLink}{\operatorname{Lk}_{\hor}}
  \newvariable{\FaceLink}{\operatorname{Lk}_{\partial}}
  \newvariable{\CofaceLink}{\operatorname{Lk}_{\delta}}

  \newvariable{\TheHemisphereCx}{H}

  \newvariable{\Homotopy}{\pi}
  \newvariable{\TheDimension}{n}

  \newvariable{\TheIndex}{i}
  \newvariable{\AltIndex}{j}
  \newvariable{\TheLastIndex}{k}
  \newvariable{\AltLastIndex}{l}

  \newvariable{\TheConst}{N}

  \newvariable{\TheNumber}{n}

  \newcommand{\CAT}[1]{\ensuremath{\operatorname{CAT}(#1)}}

  \newvariable{\OrthProj}{\operatorname{proj}}

  \newcommand{\BarycenterOf}[1]{\topcirc{#1}}
  \newvariable{\Boundary}{\partial}

  \newcommand{\uprel}{\nearrow}
  \newcommand{\downrel}{\searrow}
  \newcommand{\moveto}{\rightarrow}

  \newvariable{\ThePEC}{X}
  \newvariable{\ThePecCell}{c}
  \newvariable{\ThePecVertex}{v}
  \newvariable{\TheMorseFct}{h}
  \newvariable{\AltMorseFct}{\tilde{h}}
  \newvariable{\TheLevel}{r}
  \newvariable{\AltLevel}{s}

  \newvariable{\PiHalf}{\frac{\pi}{\Two}}
  \newvariable{\TheSphBuild}{\Delta}
  \newvariable{\TheSphChamber}{C}
  \newvariable{\TheSphSimplex}{\tau}
  \newvariable{\TheSphVertex}{v}
  \newvariable{\TheSphPoint}{x}
  \newvariable{\AltSphPoint}{y}
  \newvariable{\Distance}{d}
  \newvariable{\TheNumFactors}{k}
  \newvariable{\hor}{\mathrm{hor}}
  \newvariable{\ver}{\mathrm{ver}}
  \newcommand{\join}{*}
  \newvariable{\TheDim}{\operatorname{dim}}

  \newvariable{\TheReals}{\RRR}
  \newvariable{\TheIntegers}{\ZZZ}
  \newvariable{\TheEpsilon}{\varepsilon}
  \newvariable{\TheAffBuild}{X}
  \newvariable{\TheAffApp}{\Sigma}
  \newvariable{\AltAffApp}{\Sigma}
  \newvariable{\TheAffChamber}{C}
  \newvariable{\TheAffPoint}{x}
  \newvariable{\TheAffVertex}{v}
  \newvariable{\AltAffVertex}{w}
  \newvariable{\TheVertexSet}{V}
  \newvariable{\TheHeight}{\beta}
  \newvariable{\TheDepth}{\operatorname{dp}}
  \newvariable{\TheInfChamber}{C_\infty}
  \newvariable{\TheInfPoint}{e_\infty}
  \newvariable{\TheSimplex}{\tau}
  \newvariable{\AltSimplex}{\sigma}
  \newvariable{\TheFace}{\sigma}
  \newvariable{\TheCoface}{\xi}
  \newvariable{\TheCell}{\tau}

  \newvariable{\TheSpine}{Y}
  \newvariable{\TheCollapse}{Z}
  \newvariable{\TheHorosphere}{S}
  \newvariable{\TheHoroball}{H}
  \newvariable{\TheHoroballSet}{\mathcal{H}}

  \newvariable{\TheDivBuild}{\topcirc{\TheAffBuild}}
  \newvariable{\TheDivVertex}{\topcirc{\TheSimplex}}
  \newvariable{\AltDivVertex}{\topcirc{\AltSimplex}}

  \newcommand{\faceof}{\leq}
  \newcommand{\strictfaceof}{<}

  \newvariable{\TheGradient}{\nabla\TheHeight}

  \newvariable{\TheCoefficient}{a}

  \newvariable{\TheImage}{O}
  \newvariable{\TheOrderEmbedding}{\kappa}
  \newvariable{\TheBelow}{b}

  \newvariable{\TheDisk}{B}

\begin{document}
  \title{Connectivity Properties of Horospheres in
         Euclidean Buildings and Applications to Finiteness
         Properties of Discrete Groups}
  \author{Kai-Uwe Bux \and Kevin Wortman\thanks{
    The second author was partially supported by NSF
    grant DMS~-~0750032%
  }}
  \date{\datum{\PaperDate}}
  \maketitle

  \begin{abstract}
    Let $\TheGroupOf{\OkaRing}$ be an $\ThePrimeSet$-arithmetic subgroup
    of a connected, absolutely
    almost simple linear algebraic group $\TheGroup$
    over a global function field $\TheField$.
    We show that the sum of local ranks of $\TheGroup$
    determines the homological finiteness
    properties of $\TheGroupOf{\OkaRing}$ provided the
    $\TheField$-rank of $\TheGroup$ is $\One$. This shows that
    the general upper bound for the finiteness length of
    $\TheGroupOf{\OkaRing}$ established in an earlier paper
    is sharp in this case.

    The geometric analysis underlying our result determines the
    conectivity properties of horospheres in thick
    Euclidean buildings.
  \end{abstract}

  \section{Introduction}
    \label{sec:introduction}
    Let $\TheField$ be a global function field and suppose
    $\TheGroup$ is a connected, noncommutative, absolutely
    almost simple $\TheField$-group. Let $\ThePrimeSet$
    be a finite, nonempty set of pairwise inequivalent valuations
    on $\TheField$. We let $\OkaRing$ be the ring of
    $\ThePrimeSet$-integers in $\TheField$. We denote
    the completion of $\TheField$ with respect to
    $\ThePrime\in\ThePrimeSet$ by
    $\TheField[\ThePrime]$.
    We let
    \(
      \SumOfLocalRanksFctOf{\TheGroup,\ThePrimeSet}
      =
      \Sum[\ThePrime\in\ThePrimeSet]{
        \RankOf[{\TheField[\ThePrime]}]{\TheGroup}
      }
      .
    \)

    In \cite{Bux.Wortman:2007}, we proved:
    \begin{theorem}\label{thm:upper_bound}
      Suppose $\TheField$, $\ThePrimeSet$, and $\TheGroup$ are
      as above. If $\RankOf[\TheField]{\TheGroup}>\Zero$,
      then $\TheGroupOf{\OkaRing}$ is not of type
      \FType{\SumOfLocalRanksFctOf{\TheGroup,\ThePrimeSet}}.
    \end{theorem}

    Recall that a group $\AbstractGroup$ is of type
    \FType{\TheType} if there exists
    an Eilenberg-Mac\,Lane complex $\EilenbergMaclaneOf{\AbstractGroup}$
    with finite $\TheType$-skeleton. It follows that such a
    group has finitely generated homology and cohomology
    groups up to dimension $\TheType$.

    Theorem~\ref{thm:upper_bound} leads to
    the natural question of
    whether the groups from Theorem~\ref{thm:upper_bound}
    are of type
    \FType{\SumOfLocalRanksFctOf{\TheGroup,\ThePrimeSet}-\One}.
    Several results show that for special cases, the answer is yes.
    For example, Stuhler \cite{Stuhler:1980} proved that the answer
    is yes for groups of the form $\SlOf[\Two]{\OkaRing}$. And
    independent work of Abels \cite{Abels:1991} and
    Abramenko \cite{Abramenko:1996}
    has shown that the answer is yes for some higher
    rank examples. In particular, the answer is yes if $\TheGroup$
    is  a classical $\TheField[\ThePrimePower]$-group
    and
    \(
      \OkaRing=\FiniteFieldAd[\ThePrimePower]{\TheIndeterminate}
    \)
    where $\ThePrimePower$ is large depending on the rank
    of $\TheGroup$.

    In this paper, we add to the evidence above by proving:
    \begin{theorem}\label{thm:sharp}
      Suppose $\TheField$, $\ThePrimeSet$, and $\TheGroup$ are
      as in Theorem~\ref{thm:upper_bound}. If
      \(
        \RankOf[\TheField]{\TheGroup}=\One
        ,
      \)
      then $\TheGroupOf{\OkaRing}$ is of type
      \FType{\SumOfLocalRanksFctOf{\TheGroup,\ThePrimeSet}-\One}.
    \end{theorem}

    Thus, Theorem~\ref{thm:sharp} is a
    generalization of the result
    of Stuhler mentioned above. And together with the higher rank
    examples including those of Abels and Abramenko, it seems reasonable
    to make the following:
    \begin{conj}
      Suppose $\TheField$, $\ThePrimeSet$, and $\TheGroup$ are
      as in Theorem~\ref{thm:upper_bound}. If
      \(
        \RankOf[\TheField]{\TheGroup}>\Zero
        ,
      \)
      then $\TheGroupOf{\OkaRing}$ is of type
      \FType{\SumOfLocalRanksFctOf{\TheGroup,\ThePrimeSet}-\One}.
    \end{conj}

    \subsection{Background}
      See the introduction of \cite{Bux.Wortman:2007} for more
      on the background of this problem.

    \subsection{Outline of the proof of Theorem~\ref{thm:sharp}}
      The group $\TheGroupOf{\OkaRing}$ acts on the Euclidean
      building for
      \(
        \Product[\ThePrime\in\ThePrimeSet]{
          \TheGroupOf{\TheField[\ThePrime]}
        },
      \)
      which we denote by $\TheAffBuild$. Since
      $\RankOf[\TheField]{\TheGroup}=\One>\Zero$, the group
      $\TheGroupOf{\OkaRing}$ does not act cocompactly on
      $\TheAffBuild$. Nor does $\TheGroupOf{\OkaRing}$, nor
      any finite index subgroup of $\TheGroupOf{\OkaRing}$,
      act freely on $\TheAffBuild$, although $\TheGroupOf{\OkaRing}$
      does act on $\TheAffBuild$ with finite cell stabilizers.

      We apply reduction theory to obtain a subspace
      $\TheAffBuild[\Zero]\subseteq\TheAffBuild$
      on which $\TheGroupOf{\OkaRing}$ acts with compact quotient.
      The space
      $\TheAffBuild[\Zero]$ is obtained from $\TheAffBuild$ by
      removing an infinite family of pairwise disjoint horoballs.

      We use piecewise linear Morse theory to prove that
      horospheres -- the boundaries of horoballs -- appearing
      in the construction of $\TheAffBuild[\Zero]$ are
      $\SumOfLocalRanksFctOf{\TheGroup,\ThePrimeSet}-\Two$ connected.
      It follows that
      \(
        \TheAffBuild[\Zero]
      \)
      is $\SumOfLocalRanksFctOf{\TheGroup,\ThePrimeSet}-\Two$ connected,
      so $\TheGroupOf{\OkaRing}$ is of type
      \FType{\SumOfLocalRanksFctOf{\TheGroup,\ThePrimeSet}-\One}.

      Our proof that the horospheres used to define $\TheAffBuild[\Zero]$
      are $\SumOfLocalRanksFctOf{\TheGroup,\ThePrimeSet}-\Two$ connected
      makes essential use of a result of Schulz \cite{Schulz:2005} that
      analyzes connectivity properties of certain subsets of links in
      $\TheAffBuild$.

      We remark that unlike previous work on the positive
      direction for finiteness properties of arithmetic
      groups over function fields, our argument does not need assumptions
      on the geometry of the building $\TheAffBuild$ neither with
      regard to type, nor degree of thickness, nor the dimensions
      of irreducible factors.

    \subsection{Acknowledgements}
      The first named author thanks Bernd~Schulz for a thorough
      explaination of his PhD-thesis
      and Andrei~Rapinchuck for helpful conversations.

  \section{Reduction Theory}
    \label{sec:reduction_theory}
    Let $\TheField$ be a global function field and suppose
    $\TheGroup$ is a connected, noncommutative, absolutely almost
    simple $\TheField$-group. We assume that
    \(
      \RankOf[\TheField]{\TheGroup}
      =
      \One,
    \)
    we let $\TheTorus\subgroup\TheGroup$ be a maximal
    $\TheField$-split torus, and we choose a minimal
    $\TheField$-parabolic subgroup $\TheParabolic\subgroup\TheGroup$
    containing $\TheTorus$.

    Note that the group of $\TheField$-characters of $\TheParabolic$
    is infinite cyclic.
    Let $\TheCharacter$ be the generator for this group that
    is a positive
    multiple of the simple root associated with $\TheParabolic$ in
    the root system of $\TheGroup$ with respect to $\TheTorus$.

    Let $\ThePrimeSet$ be a finite, nonempty set of pairwise
    inequivalent valuations on $\TheField$.
    Any $\ThePrime\in\ThePrimeSet$
    gives a norm $\NormOf[\ThePrime]{\DummyArg}$ on $\TheField$, and
    we let $\TheField[\ThePrime]$ be the completion of $\TheField$
    with respect to this norm.  For any $\TheField$-group $\AltGroup$,
    we put
    \(
      \AltGroup[\ThePrimeSet]
      :=
      \Product[\ThePrime\in\ThePrimeSet]{
        \AltGroupOf{\TheField[\ThePrime]}
      }.
    \)

    The statement of the next result below requires the
    definition of two sets. The first is the group
    \[
      \TheParabolic[\ThePrimeSet][\Zero]
      :=
      \SetOf[{
        \FamOf[\ThePrime\in\ThePrimeSet]{\TheGroupElement[\ThePrime]}
        \in
        \TheParabolic[\ThePrimeSet]
      }]{
        \Product[\ThePrime\in\ThePrimeSet]{
          \NormOf[\ThePrime]{
            \TheCharacterOf{\TheGroupElement[\ThePrime]}
          }
        }
        =
        \One
      }.
    \]
    Second, for any $\TheCutOff>\Zero$, we put:
    \[
      \TheTorusOf{\TheCutOff}
      :=
      \SetOf[{
        \FamOf[\ThePrime\in\ThePrimeSet]{
          \TheGroupElement[\ThePrime]
        }
        \in
        \TheTorus[\ThePrimeSet]
      }]{
        \Product[\ThePrime\in\ThePrimeSet]{
          \NormOf[\ThePrime]{
            \TheCharacterOf{\TheGroupElement[\ThePrime]}
          }
        }
        \geq
        \TheCutOff
      }.
    \]
    We let $\OkaRing$ be the ring of $\ThePrimeSet$-integers in
    $\TheField$ and we recall that $\TheGroupOf{\OkaRing}$ is
    a discrete subgroup of $\TheGroup[\ThePrimeSet]$ via the
    diagonal embedding. The following
    theorem is a  well known result from reduction theory.
    \begin{theorem}\label{thm:finite_set}
      There is a finite set of representatives
      \(
        \TheFiniteSet\subset\TheGroupOf{\TheField}
      \)
      for the double coset space
      \(
        \TheGroupOf{\OkaRing}\lmod \TheGroupOf{\TheField} \rmod
        \TheParabolicOf{\TheField}
        .
      \)
      Furthermore, for any such set $\TheFiniteSet$, there is
      some number $\TheCutOff>\Zero$ and some compact set
      \(
        \TheCompactSet\subseteq\TheGroup[\ThePrimeSet]
      \)
      such that
      \[
        \TheGroup[\ThePrimeSet]
        =
        \TheGroupOf{\OkaRing}
        \TheFiniteSet
        \TheParabolic[\ThePrimeSet][\Zero]
        \TheTorusOf{\TheCutOff}
        \TheCompactSet
        .
      \]
    \end{theorem}
    \begin{proof}
      The finiteness of the double coset space
      \(
        \TheGroupOf{\OkaRing}\lmod \TheGroupOf{\TheField} \rmod
        \TheParabolicOf{\TheField}
      \)
      is the statement of \cite[Satz~8]{Behr:1969}. Behr's proof needs a
      technical hypothesis (used for \cite[Satz~5]{Behr:1969}).
      However, Harder has removed the need for that hypothesis:
      \cite[Korollar~2.2.7]{Harder:1969} can be used as a replacement
      for \cite[Satz~5]{Behr:1969} in the proof.

      The formula
      \(
        \TheGroup[\ThePrimeSet]
        =
        \TheGroupOf{\OkaRing}
        \TheFiniteSet
        \TheParabolic[\ThePrimeSet][\Zero]
        \TheTorusOf{\TheCutOff}
        \TheCompactSet
      \)
      follows from \cite[Satz~9]{Behr:1969}. Again, Behr uses a
      technical hypothesis, which has been subsequently removed.
      In this case, in addition to Harder's
      version of Behr's
      Satz~5, one needs to use a version of Behr's Satz~6 free from
      technical assumptions. This version is
      the main result of \cite{Springer:1994}. Using these
      replacements, Behr's proofs apply.

      We also remark that the discussion in
      \cite[page~52]{Harder:1969} implicitly contains a derivation of
      Theorem~\pref{thm:finite_set} in the context of Harder's
      version of reduction theory in positive characteristic.
    \end{proof}
    In the following section, $\TheFiniteSet$, $\TheCutOff$,
    and $\TheCompactSet$ are fixed and satisfy
    Theorem~\ref{thm:finite_set}.

  \section{Horoballs}
    \label{sec:horoballs}
    We denote the Euclidean building corresponding to
    $\TheGroupOf{\TheField[\ThePrime]}$ by $\TheAffBuild[\ThePrime]$
    and we let
    \(
      \TheAffBuild =
      \Product[\ThePrime\in\ThePrimeSet]{
        \TheAffBuild[\ThePrime]
      }
      .
    \)
    We also fix a vertex $\TheBaseVertex\in\TheAffBuild$.

    Let
    \(
      \TheRay \mapcolon \TheReals[\geq\Zero]
      \rightarrow
      \TheAffBuild
    \)
    be the geodesic ray with $\TheRayOf{\Zero}=\TheBaseVertex$
    and such that $\TheRayOf{\infty}$ is the center of mass
    of the cell corresponding to $\TheParabolic[\ThePrimeSet]$ in the
    Tits boundary of $\TheAffBuild$. Recall that
    this cell is the spherical join of the cells in the
    boundaries of the factors $\TheAffBuild[\ThePrime]$ corresponding to
    $\TheParabolicOf{\TheField[\ThePrime]}$.

    Let
    \(
      \TheHeight[\TheRay] \mapcolon
      \TheAffBuild \rightarrow \TheReals
    \)
    be the Busemann function induced by $\TheRay$
    normalized so that
    \(
      \TheHeightOf{\TheBaseVertex} = \Zero
      .
    \)
    Sets of the form
    \(
      \TheHeightOf[\TheRay][-\One]{
        \TheReals[\geq\TheParameter]
      }
    \)
    are called \notion{horoballs} -- or even horoballs
    based at $\TheParabolic[\ThePrimeSet]$ -- in analogy with
    symmetric spaces.
    The next couple of lemmas explain how Theorem~\ref{thm:finite_set}
    naturally identifies a collection of horoballs that cover
    $\TheAffBuild$.

    \begin{lemma}\label{lemma:five}
      For some $\TheParameter\in\TheReals$, we have
      \[
        \TheParabolic[\ThePrimeSet][\Zero]
        \TheTorusOf{\TheCutOff}
        \TheCompactSet
        \acts
        \TheBaseVertex
        \subseteq
        \TheHeightOf[\TheRay][-\One]{
          \TheReals[\geq\TheParameter]
        }.
      \]
    \end{lemma}
    \begin{proof}
      Let
      \(
        \ThePoint\in
        \TheParabolic[\ThePrimeSet][\Zero]
        \TheTorusOf{\TheCutOff}
        \TheCompactSet
        \acts
        \TheBaseVertex
        .
      \)
      Note that
      \(
        \TheHeightOf[\TheRay]{\TheBaseVertex}=\Zero
      \)
      and since $\TheCompactSet$ is compact,
      \(
        \TheCompactSet\acts\TheBaseVertex
        \subseteq
        \TheHeightOf[\TheRay][-\One]{
          \ClosedInterval{-\TheParameter[\Zero]}{\TheParameter[\Zero]}
        }
      \)
      for some $\TheParameter[\Zero]>\Zero$.
      Furthermore,
      \(
        \TheParabolic[\ThePrimeSet][\Zero]
        \TheTorusOf{\TheCutOff}
        \subseteq
        \TheParabolic[\ThePrimeSet]
      \)
      fixes $\TheRayOf{\infty}$ and it follows from the
      definition of $\TheParabolic[\ThePrimeSet][\Zero]$ and
      the fact that $\TheCharacter$ is a positive multiple of the
      simple root associated with $\TheParabolic$ that
      $\TheParabolic[\ThePrimeSet][\Zero]$ stabilizes
      \notion{horospheres based at $\TheParabolic[\ThePrimeSet]$},
      i.e., sets of the form
      \(
        \TheHeightOf[\TheRay][-\One]{\TheHoroRadius}
        .
      \)
      Similarly, the definition of $\TheTorusOf{\TheCutOff}$
      implies that there is a constant $\TheThickness[\TheCutOff]\geq\Zero$
      such that for any
      \(
        \TheGroupElement\in\TheTorusOf{\TheCutOff}
        ,
      \)
      we have
      \(
        \TheGroupElement
        \acts
        \TheHeightOf[\TheRay][-\One]{\TheHoroRadius}
        \subseteq
        \TheHeightOf[\TheRay][-\One]{
          \RightOpenInterval{
            \TheHoroRadius-\TheThickness[\TheCutOff]
          }{
            \infty
          }
        }
        .
      \)
    \end{proof}

    \begin{lemma}\label{lemma:six}
      There is a $\TheCoverDepth\in\TheReals$ such that
      \[
        \TheGroupOf{\OkaRing}\TheFiniteSet
        \acts
        \TheHeightOf[\TheRay][-\One]{
          \TheReals[\geq\TheCoverDepth]
        }
        =
        \TheAffBuild
        .
      \]
    \end{lemma}
    \begin{proof}
      Combining Theorem~\ref{thm:finite_set} with Lemma~\ref{lemma:five}
      shows that
      \[
        \TheGroup[\ThePrimeSet]
        \acts
        \TheBaseVertex
        \subseteq
        \TheGroupOf{\OkaRing}
        \TheFiniteSet
        \acts
        \TheHeightOf[\TheRay][-\One]{
          \TheReals[\geq\TheParameter]
        }
      \]
      for some $\TheParameter\in\TheReals$.
      The claim follows for some
      $\TheCoverDepth\leq\TheParameter$ since any point in
      $\TheAffBuild$ is a uniform bounded distance from
      a point in the orbit
      $\TheGroup[\ThePrimeSet]\acts\TheBaseVertex$.
    \end{proof}

    We have identified a cover of $\TheAffBuild$ by horoballs.
    That is not a very interesting fact on its own, but we will use
    it to help us prove our ultimate goal in this section,
    which is to identify a pairwise disjoint collection of horoballs
    in $\TheAffBuild$ with $\TheGroupOf{\OkaRing}$-invariant,
    cocompact complement.
    These pairwse disjoint horoballs will be retracts of the horoballs
    identified in Lemma~\ref{lemma:six}, so we will want to know that
    any such horoball with its retract horoball removed has a
    precompact image in
    \(
      \TheGroupOf{\OkaRing}\lmod\TheAffBuild
      .
    \)
    That is the goal of the corollary of the following:
    \begin{lemma}\label{lemma:seven}
      For $\TheHoroRadius\in\TheReals$,
      any finite index subgroup of the discrete group
      \(
        \TheParabolic[\ThePrimeSet][\Zero]
        \intersect
        \TheGroupOf{\OkaRing}
      \)
      acts cocompactly on the horosphere
      $\TheHeightOf[\TheRay][-\One]{\TheHoroRadius}$.
    \end{lemma}
    \begin{proof}
      $\TheParabolic[\ThePrimeSet][\Zero]$ stabilizes
      the horosphere
      $\TheHeightOf[\TheRay][-\One]{\TheHoroRadius}$
      as in the proof of Lemma~\ref{lemma:five}.

      Let $\AltPoint\in\TheHeightOf[\TheRay][-\One]{\TheHoroRadius}$
      and let $\TheAffApp[\ThePrime]\subseteq\TheAffBuild[\ThePrime]$
      be an apartment corresponding to a maximal
      $\TheField[\ThePrime]$-split torus
      $\TheTorus[\ThePrime]$
      in $\TheParabolic$
      containing $\TheTorus$.
      Let
      $\TheRadical$ be the unipotent radical of $\TheParabolic$.
      Recall that $\TheRadicalOf{\TheField[\ThePrime]}$ acts
      transitively on the set of apartments in $\TheAffBuild[\ThePrime]$
      whose boundary sphere contains the chamber corresponding
      to $\TheParabolicOf{\TheField[\ThePrime]}$.
      Since
      \(
        \TheRadicalOf{\TheField[\ThePrime]}\subgroup
        \TheParabolic[\ThePrimeSet][\Zero],
      \)
      there is some $\TheGroupElement\in\TheParabolic[\ThePrimeSet][\Zero]$
      with
      \(
        \TheGroupElement
        \acts
        \AltPoint
        \in
        \Parentheses{
          \Product[\ThePrime\in\ThePrimeSet]{
            \TheAffApp[\ThePrime]
          }
        }
        \intersect
        \TheHeightOf[\TheRay][-\One]{\TheHoroRadius}
        .
      \)

      Note that
      \(
        \Parentheses{
          \Product[\ThePrime\in\ThePrimeSet]{
            \TheAffApp[\ThePrime]
          }
        }
        \intersect
        \TheHeightOf[\TheRay][-\One]{\TheHoroRadius}
      \)
      is a codimension $\One$-subspace of
      \(
          \Product[\ThePrime\in\ThePrimeSet]{
            \TheAffApp[\ThePrime]
          }
          ,
      \)
      and that
      \(
        \Parentheses{
          \Product[\ThePrime\in\ThePrimeSet]{
            \TheTorusOf[\ThePrime]{\TheField[\ThePrime]}
          }
        }
        \intersect
        \TheParabolic[\ThePrimeSet][\Zero]
      \)
      acts cocompactly on
      \(
        \Parentheses{
          \Product[\ThePrime\in\ThePrimeSet]{
            \TheAffApp[\ThePrime]
          }
        }
        \intersect
        \TheHeightOf[\TheRay][-\One]{\TheHoroRadius}
        .
      \)
      We have shown that
      $\TheParabolic[\ThePrimeSet][\Zero]$ acts cocompactly
      on the horosphere
      \(
        \TheHeightOf[\TheRay][-\One]{\TheHoroRadius}
        .
      \)

      Now, we observe that
      \(
        \TheParabolic[\ThePrimeSet][\Zero]
        \intersect
        \TheGroupOf{\OkaRing}
      \)
      is cocompact in $\TheParabolic[\ThePrimeSet][\Zero]$:
      the parabolic group $\TheParabolic[\ThePrimeSet][\Zero]$
      decomposes as a product
      $\TheRadical[\ThePrimeSet]
      \TheLevi[\ThePrimeSet]
      \TheTorus[\ThePrimeSet][\Zero]$ where
      $\TheLevi$ is
      a reductive group of $\TheField$-rank $\Zero$.
      Note that the subgroup of
      $\ThePrimeSet$-integer points in $\TheRadical[\ThePrimeSet]$
      is cocompact since $\TheRadical$ is unipotent;
      for the subgroup of $\ThePrimeSet$-integer
      points in $\TheLevi[\ThePrimeSet]$,
      cocompactness follows
      from \cite[Korollar~2.2.7]{Harder:1969}; and for
      $\TheTorus[\ThePrimeSet][\Zero]$,
      cocompactness of the $\ThePrimeSet$-integer
      subgroup follows from Dirichlet's unit theorem.
      Therefore,
      \(
        \TheParabolic[\ThePrimeSet][\Zero]
        \intersect
        \TheGroupOf{\OkaRing}
      \)
      is cocompact in $\TheParabolic[\ThePrimeSet][\Zero]$
      and thus so is any finite index subgroup
      of
      \(
        \TheParabolic[\ThePrimeSet][\Zero]
        \intersect
        \TheGroupOf{\OkaRing}
        .
      \)
      The claim now follows.
    \end{proof}
    \begin{cor}\label{cor:eight}
      For any $\TheFiniteSetElement\in\TheFiniteSet$
      and any $\TheSliceSize\geq\Zero$, the quotient
      \[
        \lquot{
          \TheGroupOf{\OkaRing}
        }{
          \TheFiniteSetElement
          \acts
          \TheHeightOf[\TheRay][-\One]{
            \ClosedInterval{-\TheSliceSize}{\TheSliceSize}
          }
        }
        \subseteq
        \lquot{
          \TheGroupOf{\OkaRing}
        }{
          \TheAffBuild
        }
      \]
      is compact.
    \end{cor}
    \begin{proof}
      For $\TheFiniteSetElement=\One$, the
      claim is immediate from Lemma~\ref{lemma:seven}.
      If $\TheFiniteSetElement\neq\One$, we replace the role of
      \(
        \TheParabolic[\ThePrimeSet][\Zero]\intersect
        \TheGroupOf{\OkaRing}
        =
        \TheParabolicOf{\OkaRing}
      \)
      in Lemma~\ref{lemma:seven} with
      \(
        \TheFiniteSetElement
        \TheParabolicOf{\OkaRing}
        \TheFiniteSetElement[][-\One]
        \intersect
        \TheGroupOf{\OkaRing}
        .
      \)
      This is a finite index subgroup of
      \(
        \TheFiniteSetElement
        \TheParabolicOf{\OkaRing}
        \TheFiniteSetElement[][-\One]
      \)
      since
      $\TheFiniteSetElement\in\TheGroupOf{\TheField}$,
      see e.g.\ \cite[Lemma~3.1.1(iv)]{Margulis:1991}.
    \end{proof}

    We let
    \(
      \TheHeight[\TheFiniteSetElement\acts\TheRay]
      \mapcolon
      \TheAffBuild
      \rightarrow
      \TheReals
    \)
    be the Busemann function for the geodesic ray
    \(
      \TheFiniteSetElement\acts\TheRay.
    \)
    Thus,
    \(
      \TheFiniteSetElement\acts
      \TheHeightOf[\TheRay][-\One]{\TheHoroRadius}
      =
      \TheHeightOf[\TheFiniteSetElement\acts\TheRay][-\One]{
        \TheHoroRadius
      }.
    \)
    \begin{lemma}\label{lemma:nine}
      There is some
      \(
        \TheHoroDepth[\TheFiniteSetElement]
        >
        \Zero
      \)
      such that
      \(
        \TheGroupOf{\OkaRing}\acts\TheBaseVertex
        \intersect
        \TheHeightOf[\TheFiniteSetElement\acts\TheRay][-\One]{
          \RightOpenInterval{
            \TheHoroDepth[\TheFiniteSetElement]
          }{
            \infty
          }
        }
        =
        \emptyset
        .
      \)
    \end{lemma}
    \begin{proof}
      Let
      \(
        \TheLattice[\TheFiniteSetElement]
        =
        \TheFiniteSetElement
        \TheParabolicOf{\OkaRing}
        \TheFiniteSetElement[][-\One]
        \intersect
        \TheGroupOf{\OkaRing}
      \)
      and let
      $\TheDistance[\TheRadius]$ be the Hausdorff distance between
      the orbit
      \(
        \TheLattice[\TheFiniteSetElement]
        \TheFiniteSetElement
        \acts
        \TheRayOf{\TheRadius}
      \)
      and the horosphere
      \(
        \TheHeightOf[\TheFiniteSetElement\acts\TheRay][-\One]{
          \TheRadius
        }
        .
      \)
      Note that
      \(
        \TheDistance[{\TheRadius[\Two]}]
        \leq
        \TheDistance[{\TheRadius[\One]}]
      \)
      when
      \(
        \TheRadius[\One] < \TheRadius[\Two]
        .
      \)

      Let $\TheTorusElement\in\TheTorus[\ThePrimeSet]$ be defined
      by $\TheTorusElement=\FamOf[\ThePrime\in\ThePrimeSet]{\TheTorusElement[\ThePrime]}$
      where $\TheTorusElement[\ThePrime]\in\TheTorusOf{\TheField[\ThePrime]}$
      is such that
      $\NormOf[\ThePrime]{\TheCharacterOf{\TheTorusElement[\ThePrime]}}>\One$,
      and let
      $\TheTorusElement[\TheFiniteSetElement]
      =
      \TheFiniteSetElement\TheTorusElement\TheFiniteSetElement[][-\One]$.
      Since
      \(
        \TheFiniteSetElement\TheTorusElement[\ThePrime]
        \TheFiniteSetElement[][-\One]
      \)
      acts by translations on
      $\TheFiniteSetElement\acts\TheAffApp[\ThePrime]$,
      we have
      \(
        \TheHeightOf[\TheFiniteSetElement\acts\TheRay]{
          \TheTorusElement[\TheFiniteSetElement][\TheExponent]
          \acts
          \TheBaseVertex
        }
        =
        \TheExponent\TheLag
        +
        \TheHeightOf[\TheFiniteSetElement\acts\TheRay]{
          \TheBaseVertex
        }
      \)
      for some $\TheLag>\Zero$.
      Note that for any
      \(
        \TheConjugateUnip\in
        \TheFiniteSetElement
        \TheRadicalOf{\OkaRing}
        \TheFiniteSetElement[][-\One]
        \intersect
        \TheGroupOf{\OkaRing}
        ,
      \)
      the sequence
      \(
        \TheTorusElement[\TheFiniteSetElement][-\TheExponent]
        \TheConjugateUnip
        \TheTorusElement[\TheFiniteSetElement][\TheExponent]
      \)
      converges to $\TheIdentityElement$ in
      $\TheGroup[\ThePrimeSet]$. By
      \cite[Theorem~I.1.12]{Raghunathan:1972}, the sequence
      \(
        \TheGroupOf{\OkaRing}
        \TheTorusElement[\TheFiniteSetElement][-\TheExponent]
        \TheConjugateUnip
        \TheTorusElement[\TheFiniteSetElement][\TheExponent]
      \)
      contains no convergent subsequence in
      \(
        \TheGroupOf{\OkaRing}
        \lmod
        \TheGroup[\ThePrimeSet]
        .
      \)
      It follows that there is some $\TheBigExponent$ such
      that
      \(
        \DistanceOf{
          \TheTorusElement[\TheFiniteSetElement][\TheExponent]
          \acts
          \TheBaseVertex
          ,
          \TheGroupOf{\OkaRing}
          \acts
          \TheBaseVertex
        }
        >
        \TheLag + \TheDistance[\Zero]
      \)
      for all $\TheExponent\geq\TheBigExponent$.

      Let
      \(
        \TheHoroDepth[\TheFiniteSetElement]
        =
        \MaxOf{
          \TheBigExponent\TheLag,
          \TheHeightOf[\TheFiniteSetElement\acts\TheRay]{
            \TheBaseVertex
          }
        }
      \)
      and suppose
      \(
        \ThePoint
        \in
        \TheHeightOf[\TheFiniteSetElement\acts\TheRay][-\One]{
          \RightOpenInterval{\TheHoroDepth[\TheFiniteSetElement]}{\infty}
        }
        .
      \)
      We claim that
      \(
        \ThePoint\not\in\TheGroupOf{\OkaRing}\acts\TheBaseVertex
        .
      \)
      Indeed, we can choose $\TheExponent\geq\TheBigExponent$ such
      that
      \[
        \AbsValueOf{
          \TheHeightOf[\TheFiniteSetElement\acts\TheRay]{
            \TheTorusElement[\TheFiniteSetElement][\TheExponent]
            \acts
            \TheBaseVertex
          }
          -
          \TheHeightOf[\TheFiniteSetElement\acts\TheRay]{
            \ThePoint
          }
        }
        <
        \TheLag
      \]
      so we can choose
      \(
        \TheLatticeElement\in\TheLattice[\TheFiniteSetElement]
      \)
      such that
      \[
        \DistanceOf{
          \TheTorusElement[\TheFiniteSetElement][\TheExponent]
          \acts
          \TheBaseVertex
          ,
          \TheLatticeElement
          \acts
          \ThePoint
        }
        \leq
        \TheLag+\TheDistance[\Zero]
        .
      \]
      Thus
      \(
        \TheLatticeElement
        \acts
        \ThePoint
        \not\in
        \TheGroupOf{\OkaRing}\acts
        \TheBaseVertex
      \)
      and so
      \(
        \ThePoint\not\in
        \TheGroupOf{\OkaRing}\acts
        \TheBaseVertex
        .
      \)
    \end{proof}

    We let $\ThePush>\Zero$ be the maximum of the distances from
    points in the horosphere
    $\TheHeightOf[\TheFiniteSetElement\acts\TheRay][-\One]{
      \TheHoroDepth[\TheFiniteSetElement]
    }$
    to the orbit
    $\TheGroupOf{\OkaRing}\acts\TheBaseVertex$
    as $\TheFiniteSetElement$ ranges through $\TheFiniteSet$.
    Note that $\ThePush$ is finite by Lemma~\ref{lemma:seven}.
    \begin{lemma}\label{lemma:ten}
      Let $\TheArithmeticElement\in\TheGroupOf{\OkaRing}$ and
      suppose
      \[
        \TheArithmeticElement
        \acts
        \TheHeightOf[\TheFiniteSetElement\acts\TheRay][-\One]{
          \RightOpenInterval{
            \TheHoroDepth[\TheFiniteSetElement]
            +
            \ThePush
          }{
            \infty
          }
        }
        \intersect
        \TheHeightOf[\AltFiniteSetElement\acts\TheRay][-\One]{
          \RightOpenInterval{
            \squish{\TheHoroDepth[\AltFiniteSetElement]}
            +
            \ThePush
          }{
            \infty
          }
        }
        \neq
        \emptyset
        .
      \]
      Then
      \(
        \TheArithmeticElement\in
        \ParenthesesOf{
          \TheFiniteSetElement
          \TheParabolic
          \TheFiniteSetElement[][-\One]
        }{
          \OkaRing
        }
        ,
      \)
      and
      $\TheFiniteSetElement=\AltFiniteSetElement$,
      and
      \[
        \TheArithmeticElement\acts
        \TheHeightOf[\TheFiniteSetElement\acts\TheRay][-\One]{
          \RightOpenInterval{
            \TheHoroDepth[\TheFiniteSetElement]
            +
            \ThePush
          }{
            \infty
          }
        }
        =
        \TheHeightOf[\AltFiniteSetElement\acts\TheRay][-\One]{
          \RightOpenInterval{
            \squish{\TheHoroDepth[\AltFiniteSetElement]}
            +
            \ThePush
          }{
            \infty
          }
        }
        .
      \]
    \end{lemma}
    \begin{proof}
      Let
      $\AltAffApp\subseteq\TheAffBuild$
      be an apartment
      whose
      boundary sphere at infinity
      contains the cells corresponding
      to
      \(
          \TheArithmeticElement
          \TheFiniteSetElement
          \TheParabolic[\ThePrimeSet]
          \TheFiniteSetElement[][-\One]
          \TheArithmeticElement[][-\One]
      \)
      and
      \(
        \AltFiniteSetElement
        \TheParabolic[\ThePrimeSet]
        \AltFiniteSetElement[][-\One]
        .
      \)
      If
      \(
        \TheArithmeticElement
        \TheFiniteSetElement
        \TheParabolic[\ThePrimeSet]
        \TheFiniteSetElement[][-\One]
        \TheArithmeticElement[][-\One]
        \neq
        \AltFiniteSetElement
        \TheParabolic[\ThePrimeSet]
        \AltFiniteSetElement[][-\One]
      \)
      then $\RankOf[\TheField]{\TheGroup}=\One$ implies
      that these are opposite cells at infinity, and thus
      the triple intersection
      \[
        \AltAffApp
        \intersect
        \TheArithmeticElement
        \acts
        \TheHeightOf[\TheFiniteSetElement\acts\TheRay][-\One]{
          \RightOpenInterval{
            \TheHoroDepth[\TheFiniteSetElement]
            +
            \ThePush
          }{
            \infty
          }
        }
        \intersect
        \TheHeightOf[\AltFiniteSetElement\acts\TheRay][-\One]{
          \RightOpenInterval{
            \squish{\TheHoroDepth[\AltFiniteSetElement]}
            +
            \ThePush
          }{
            \infty
          }
        }
      \]
      is contained in a metric neighborhood of a
      hyperplane in $\AltAffApp$.

      We choose
      \[
        \ThePoint\in
        \AltAffApp
        \intersect
        \TheArithmeticElement
        \acts
        \TheHeightOf[\TheFiniteSetElement\acts\TheRay][-\One]{
          \RightOpenInterval{
            \TheHoroDepth[\TheFiniteSetElement]
            +
            \ThePush
          }{
            \infty
          }
        }
        \intersect
        \TheHeightOf[\AltFiniteSetElement\acts\TheRay][-\One]{
          \squish{\TheHoroDepth[\AltFiniteSetElement]}
        }
        .
      \]
      It follows from the choice of $\ThePush$ that there is some
      \(
        \AltPoint\in
        \TheGroupOf{\OkaRing}\acts\TheBaseVertex
      \)
      such that
      \(
        \DistanceOf{\ThePoint,\AltPoint}
        \leq\ThePush
        .
      \)
      Therefore,
      \(
        \TheHeightOf[
          \TheArithmeticElement\TheFiniteSetElement\acts\TheRay
        ]{
          \AltPoint
        }
        \geq
        \TheHoroDepth[\TheFiniteSetElement]
      \)
      which contradicts Lemma~\ref{lemma:nine}.
      Hence
      \(
        \TheArithmeticElement
        \TheFiniteSetElement
        \TheParabolic[\ThePrimeSet]
        \TheFiniteSetElement[][-\One]
        \TheArithmeticElement[][-\One]
        =
        \AltFiniteSetElement
        \TheParabolic[\ThePrimeSet]
        \AltFiniteSetElement[][-\One]
      \)
      which is to say that
      $\TheFiniteSetElement=\AltFiniteSetElement$ and
      \(
        \TheArithmeticElement
        \in
        \ParenthesesOf{
          \TheFiniteSetElement
          \TheParabolic[\ThePrimeSet]
          \TheFiniteSetElement[][-\One]
        }{
          \OkaRing
        }
        .
      \)

      Furthermore, $\TheArithmeticElement$ preserves distances
      from $\TheGroupOf{\OkaRing}\acts\TheBaseVertex$,
      so the result follows.
    \end{proof}

    Let
    $\TheHoroCollection=\TheGroupOf{\OkaRing}\acts
    \SetOf[{
      \TheHeightOf[\TheFiniteSetElement\acts\TheRay][-\One]{
        \OpenInterval{
          \TheHoroDepth[\TheFiniteSetElement]+\ThePush
        }{\infty}
      }
    }]{
      \TheFiniteSetElement\in\TheFiniteSet
    }.$
    This is a collection of open horoballs.
    \begin{theorem}\label{thm:geometry}
      $\TheHoroCollection$ is a
      collection of pairwise disjoint horoballs;
      $\TheAffBuild\setminus\TheHoroCollection$
      is $\TheGroupOf{\OkaRing}$-invariant and cocompact.
    \end{theorem}
    \begin{proof}
      The result follows from the definition of $\TheHoroCollection$,
      Lemma~\ref{lemma:ten}, and Corollary~\ref{cor:eight}.
    \end{proof}

  \section{Connectivity of Horospheres in General Position}%
    \label{sec:general_pos}
    Let $\TheSphBuild$ be a spherical building. We consider
    $\TheSphBuild$ as a metric space with the angular
    metric $\Distance$, i.e., every apartment is a unit sphere.
    For any point $\TheSphPoint\in\TheSphBuild$, we define
    the \notion{closed hemisphere complex}
    \(
      \TheSphBuildOf[][\geq\PiHalf]{\TheSphPoint}
    \)
    to be the subcomplex spanned by all vertices in the set
    \(
      \SetOf[\AltSphPoint\in\TheSphBuild]{
        \DistanceOf{\AltSphPoint,\TheSphPoint}
        \geq\PiHalf
      }      ,
    \)
    the \notion{open hemisphere complex}
    \(
      \TheSphBuildOf[][>\PiHalf]{\TheSphPoint}
    \)
    as the subcomplex spanned by all vertices in the set
    \(
      \SetOf[\AltSphPoint\in\TheSphBuild]{
        \DistanceOf{\AltSphPoint,\TheSphPoint}
        >\PiHalf
      }
      ,
    \)
    and we define the the \notion{equator} as the set
    of points
    \(
      \TheSphBuildOf[][=\PiHalf]{\TheSphPoint}
      :=
      \SetOf[\AltSphPoint\in\TheSphBuild]{
        \DistanceOf{\AltSphPoint,\TheSphPoint}
        =\PiHalf
      }
      .
    \)
    Recall that $\TheSphBuild$ decomposes uniquely as the
    spherical join of irreducible factors
    \[
      \TheSphBuild
      =
      \TheSphBuild[\One]\join\cdots\join\TheSphBuild[\TheNumFactors]
      .
    \]
    where the decomposition is determined by the geometry of
    chambers as follows:
    \begin{lemma}\label{factors}
      In an irreducible spherical building, every edge has angular
      length strictly less than $\PiHalf$. Consequently,
      in any spherical building, an edge has angular length
      $\PiHalf$ if and only if it joins two vertices from
      different irreducible factors.
    \end{lemma}
    \begin{proof}
      We start with the following observation from spherical
      geometry: suppose all edges and angles in a spherical triangle
      are at most $\PiHalf$; if one edge has length exactly  $\PiHalf$
      then so has at least one of the other edges (in fact, also at
      least two of the angles will be right angles).

      We apply this observation to the vertices of a chamber
      $\TheSphChamber$. Let
      \(
        \TheSphVertex[\One],\TheSphVertex[\Two],\ldots,
        \TheSphVertex[\TheLastIndex]
      \)
      be a maximal collection of vertices that have pairwise distance
      $\PiHalf$. Then every other vertex has distance strictly less
      than $\PiHalf$ to at least one of them (by maximality of the
      collection) and therefore to exactly one of them (by the observation).
      It follows that being of distance strictly less than $\PiHalf$
      is an equivalence relation on the set of vertices of $\TheSphChamber$
      with $\TheLastIndex$ equivalence classes. This defines a decomposition
      of the underlying Coxeter complex as a spherical join and induces
      a decomposition of the building into irreducible factors.
    \end{proof}

    Let $\TheSphBuildOf[\hor]{\TheSphPoint}$ be the join of all
    irreducible factors completely contained in the equator,
    and let $\TheSphBuildOf[\ver]{\TheSphPoint}$ be the join of
    the other factors. Clearly
    \(
      \TheSphBuild =
      \TheSphBuildOf[\hor]{\TheSphPoint}
      \join
      \TheSphBuildOf[\ver]{\TheSphPoint}
      .
    \)
    \begin{lemma}\label{is_equ}
      Let $\TheSphChamber\subseteq\TheSphBuild$ be a chamber.
      Then for any equatorial simplex
      \(
        \TheSphSimplex\subseteq\TheSphBuildOf[][=\PiHalf]{\TheSphPoint}
      \)
      contained in $\TheSphChamber$,
      the following are equivalent:
      \begin{enumerate}
        \item\label{equ_a}
          We have
          \(
            \TheSphSimplex\subseteq\TheSphBuild[\hor]
            .
          \)
        \item\label{equ_b}
          The simplex $\TheSphSimplex$ has Hausdorff distance
          $\PiHalf$ from any non-equatorial vertex in $\TheSphBuild$.
        \item\label{equ_c}
          The simplex $\TheSphSimplex$ has Hausdorff distance
          $\PiHalf$ from any non-equatorial vertex in
          $\TheSphChamber$.
      \end{enumerate}
    \end{lemma}
    \begin{proof}
      The implications
      \pref{equ_a}$\Rightarrow$\pref{equ_b}
      and
      \pref{equ_b}$\Rightarrow$\pref{equ_c}
      are obvious. It remains to show that
      \pref{equ_c} implies \pref{equ_a}.

      We will show that $\TheSphVertex\not\in\TheSphBuild[\ver]$ for each
      vertex $\TheSphVertex\in\TheSphSimplex$. This implies
      that all vertices of $\TheSphSimplex$ belong to $\TheSphBuild[\hor]$
      and thus proves the claim.

      Let $\TheSphVertex$ be a vertex of $\TheSphSimplex$ and
      let $\TheSphBuild[\TheIndex]$ be a vertical irreducible factor of
      $\TheSphBuild$. Note that a chamber in $\TheSphBuild[\TheIndex]$
      cannot have all its vertices in the equator. It follows that
      $\TheSphChamber\intersect\TheSphBuild[\TheIndex]$ contains
      a non-equatorial vertex. Since this vertex is connected to
      $\TheSphVertex$ by an edge of length
      $\PiHalf$, it follows from Fact~\pref{factors} that
      $\TheSphVertex\not\in\TheSphBuild[\TheIndex]$. Since the same
      argument proves that $\TheSphVertex$ is not in any irreducible
      factor of $\TheSphBuild[\ver]$, we have
      \(
        \TheSphVertex\not\in\TheSphBuild[\ver].
      \)
    \end{proof}

    Connectivity properties of hemisphere complexes are
    given by the following:
    \begin{theorem}[{Bernd~Schulz, \cite{Schulz:2005}}]\label{bernd}
      Assume that $\TheSphBuild$ is a thick spherical
      building.
      The closed hemisphere complex
      \(
        \TheSphBuildOf[][\geq\PiHalf]{\TheSphPoint}
      \)
      is $\Parentheses{\TheDimOf{\TheSphBuild}-\One}$-connected.

      The open hemisphere complex
      \(
        \TheSphBuildOf[][>\PiHalf]{\TheSphPoint}
      \)
      is $\Parentheses{\TheDimOf{\TheSphBuildOf[\ver]{\TheSphPoint}}-\One}$-connected.
    \end{theorem}

    As a first application, we shall deduce the connectivity
    of horospheres ``in general position''.
    \begin{prop}\label{general_pos}
      Let $\TheAffBuild$ be a thick Euclidean building and let
      \[
        \TheHeight \mapcolon
        \TheAffBuild \longrightarrow \TheReals
      \]
      be a Busemann function that is non-constant on each edge
      of $\TheAffBuild$. Then, any horosphere
      \(
        \TheHeightOf[][-\One]{\TheLevel}
      \)
      is $\Parentheses{\TheDimOf{\TheAffBuild}-\Two}$-connected.
    \end{prop}
    Before we embark on the proof, we need to state a version
    of the Morse-Lemma that fuels Bestvina-Brady type combinatorial
    Morse-theory as introduced in \cite{Bestvina.Brady:1997}.
    Let $\ThePEC$ be a piecewise Euclidean complex and let
    \(
      \TheMorseFct \mapcolon \ThePEC \rightarrow \TheReals
    \)
    be a function that is affine on cells and non-constant on
    edges. The \notion{descending link}
    \(
      \DescLinkOf{\ThePecVertex}
    \)
    of a vertex
    $\ThePecVertex\in\ThePEC$ is the subcomplex
    of $\LinkOf{\ThePecVertex}$ defined by all cells $\ThePecCell$
    in $\ThePEC$ containing $\ThePecVertex$ as the point where
    $\TheMorseFct$ attains its maximum on $\ThePecCell$.
    \begin{NewTh}<statement>[Satz]{Morse Lemma}
      Let $\ThePEC$ and $\TheMorseFct$ be as above, and
      let $\TheLevel<\AltLevel$ be real numbers chosen such that
      the preimage
      \(
        \TheMorseFctOf[][-\One]{
          \ClosedInterval{\TheLevel}{\AltLevel}
        }
      \)
      does not contain a complete edge. Then
      the sublevel set
      \(
        \TheMorseFctOf[][-\One]{
          \LeftOpenInterval{-\infty}{\AltLevel}
        }
      \)
      is homotopy equivalent to the sublevel set
      \(
        \TheMorseFctOf[][-\One]{
          \LeftOpenInterval{-\infty}{\TheLevel}
        }
      \)
      with descending links of vertices in
      \(
        \TheMorseFctOf[][-\One]{
          \LeftOpenInterval{\TheLevel}{\AltLevel}
        }
      \)
      conned off.
    \end{NewTh}
    We shall not give a proof of the Morse-Lemma here since we
    will prove a slightly more general version later. We just remark
    that the version above is essentially the Morse-Lemma from
    \cite[Lemma~7]{Bux.Gonzalez:1999}.
    \begin{observation}\label{parallel}
      Let $\TheInfChamber$ be a chamber of the spherical building
      at infinity that contains the end $\TheInfPoint$.
      Note that $\TheAffBuild$ is covered by apartments
      containing $\TheInfChamber$ and that $\TheHeight$ is
      affine on all those apartments.
      In each such apartment, there are only finitely many edges
      up to translation. Moreover, any two such
      apartments have a common sector representing $\TheInfChamber$.
      Thus, there are only finitely many ``parallelism classes'' of
      edges in $\TheAffBuild$.\qed
    \end{observation}
    \begin{proof}[of Proposition~\ref{general_pos}]
      Let $\TheAffVertex\in\TheAffBuild$ be a vertex.
      The link of $\TheAffVertex$ is a spherical building
      $\TheSphBuild := \LinkOf{\TheAffVertex}$.
      There is a unique
      geodesic ray issuing from $\TheAffVertex$ toward $\TheInfPoint$.
      This geodesic represents the gradient
      \(
        \TheGradient\in\TheSphBuild=\LinkOf{\TheAffVertex}
      \)
      of the Busemann function defined by $\TheInfPoint$.
      Observe that directions issuing from $\TheAffVertex$ are
      descending if they span an obtuse angle with the gradient.
      It follows that the descending link
      $\DescLinkOf{\TheAffVertex}$ is the open hemisphere complex
      \(
        \TheSphBuildOf[][>\PiHalf]{\TheGradient}
        .
      \)
      Since there are no horizontal edges, the open hemisphere
      complex and the closed hemisphere complex coincide. It
      follows that descending links are
      \(
        \Parentheses{\TheDimOf{\TheAffBuild}-\Two}
      \)-connected.

      It follows form Observation~\pref{parallel}
      that there is a constant $\TheEpsilon>\Zero$
      so that for any two vertices $\TheAffVertex,\AltAffVertex
      \in\TheAffBuild$ joined by an edge, we have
      \(
        \AbsValueOf{
          \TheHeightOf{\TheAffVertex}
          -
          \TheHeightOf{\AltAffVertex}
        }
        >
        \TheEpsilon
        .
      \)

      Our choice of $\TheEpsilon$ ensures that no preimage
      \(
        \TheHeightOf[][-\One]{
          \ClosedInterval{\AltLevel}{\AltLevel+\TheEpsilon}
        }
      \)
      contains a complete edge. Thus, the Morse-Lemma implies that,
      for any $\AltLevel\in\TheReals$,
      \(
        \TheHeightOf[][-\One]{
          \LeftOpenInterval{-\infty}{\AltLevel+\TheEpsilon}
        }
      \)
      is homotopy equivalent to
      \(
        \TheHeightOf[][-\One]{
          \LeftOpenInterval{-\infty}{\AltLevel}
        }
      \)
      with descending links conned off. As descending links
      are
      \(
        \Parentheses{\TheDimOf{\TheAffBuild}-\Two}
      \)-connected,
      we find that the inclusion
      \[
        \TheHeightOf[][-\One]{
          \LeftOpenInterval{-\infty}{\AltLevel}
        }
        \monorightarrow
        \TheHeightOf[][-\One]{
          \LeftOpenInterval{-\infty}{\AltLevel+\TheEpsilon}
        }
      \]
      induces isomorphisms in $\Homotopy[\TheDimension]$ for
      $\TheDimension\leq\TheDimOf{\TheAffBuild}-\Two$.
      Iterating, we obtain that for any $\AltLevel>\TheLevel$
      the inclusion
      \[
        \TheHeightOf[][-\One]{
          \LeftOpenInterval{-\infty}{\TheLevel}
        }
        \monorightarrow
        \TheHeightOf[][-\One]{
          \LeftOpenInterval{-\infty}{\TheLevel+\TheEpsilon}
        }
        \monorightarrow
        \TheHeightOf[][-\One]{
          \LeftOpenInterval{-\infty}{\TheLevel+\Two\TheEpsilon}
        }
        \monorightarrow\cdots\monorightarrow
        \TheAffBuild
      \]
      induces isomorphisms in $\Homotopy[\TheDimension]$ for
      $\TheDimension\leq\TheDimOf{\TheAffBuild}-\Two$.
      Since $\TheAffBuild$ is contractible, it follows that
      sublevel sets are
      \(
        \Parentheses{\TheDimOf{\TheAffBuild}-\Two}
      \)-connected.
    \end{proof}

  \section{Toward a Secondary Morse Function}
    \label{sec:secondary}
    Ultimately, we want to deal with horospheres that are not necessarily
    in general position, i.e., the corresponding Busemann function
    might be constant on some edges. To overcome this obstacle,
    we construct a secondary Morse function that will allow us
    to break ties.

    Let $\TheAffBuild$ be an irreducible Euclidean building.
    The link of any simplex $\TheSimplex$, is the union of all
    those simplices $\AltSimplex$ disjoint from $\TheSimplex$
    such that $\AltSimplex\union\TheSimplex$ is a simplex.
    This link, $\LinkOf{\TheSimplex}$ is a spherical building, and we
    may alternatively think of its points as directions issuing
    from the barycenter $\BarycenterOf{\TheSimplex}$
    of $\TheSimplex$ that are perpendicular to $\TheSimplex$ -- one
    way to make sense of perpendicularity is to recall that
    $\TheAffBuild$ is a \CAT{\Zero} space and that $\TheSimplex$
    is a convex subset. We note that this way the link
    $\LinkOf{\TheSimplex}$ is endowed with an angular metric
    so that each apartment in $\LinkOf{\TheSimplex}$ is a unit sphere.

    Let
    \(
      \TheHeight \mapcolon
      \TheAffBuild \longrightarrow \TheReals
    \)
    be a Busemann function on $\TheAffBuild$ corresponding to
    a point $\TheInfPoint$ at infinity.
    We call a simplex
    $\TheSimplex\subseteq\TheAffBuild$ \notion{horizontal} if
    $\TheHeight$ restricts to a constant map on $\TheSimplex$.
    For a horizontal simplex $\TheSimplex$, the unique geodesic ray
    from the barycenter $\BarycenterOf{\TheSimplex}$ to the
    end $\TheInfPoint$ is perpendicular to $\TheSimplex$ and thus
    determines a direction $\TheGradient\in\LinkOf{\TheSimplex}$, to
    which we refer as the \notion{gradient} of $\TheHeight$.

    For a horizontal simplex $\TheSimplex$ with link
    $\TheSphBuild:=\LinkOf{\TheSimplex}$, we define the descending
    link $\DescLinkOf[\TheHeight]{\TheSimplex}$ as the subcomplex defined by
    those simplices in $\TheAffBuild$ that contain $\TheSimplex$ as a face
    and where $\TheHeight$ is maximal exactly along the face $\TheSimplex$.
    It is obvious that
    $\DescLinkOf[\TheHeight]{\TheSimplex}$ coincides with the
    open hemisphere complex
    \(
      \TheSphBuildOf[][>\PiHalf]{\TheGradient}
      .
    \)
    We define the equatorial link
    \(
      \EquLinkOf[\TheHeight]{\TheSimplex}
      :=
      \TheSphBuildOf[][=\PiHalf]{\TheGradient}
      ,
    \)
    i.e., as the set of directions along which $\TheHeight$ does
    not change. We also
    define the vertical link as
    \(
      \VertLinkOf{\TheSimplex}
      :=
      \TheSphBuild[\ver]
    \)
    and the horizontal link as
    \(
      \HorLinkOf{\TheSimplex}
      :=
      \TheSphBuild[\hor]
      .
    \)

    \begin{lemma}\label{perpendicular}
      Let $\TheSimplex$ be a simplex in $\TheAffBuild$ and
      let $\AltSimplex[\One],\AltSimplex[\Two]$ be two
      simplices in $\LinkOf{\TheSimplex}$ that span a
      simplex
      \(
        \AltSimplex[\One]\union\AltSimplex[\Two]
        .
      \)
      Then the following are equivalent:
      \begin{enumerate}
        \item\label{perp_b}
          The simplices $\AltSimplex[\One]$ and
          $\AltSimplex[\Two]$ have distance
          $\PiHalf$ in $\LinkOf{\TheSimplex}$.
        \item\label{perp_c}
          There is a decomposition
          \(
            \LinkOf{\TheSimplex}
            =
            \TheSphBuild[\One]\join\TheSphBuild[\Two]
          \)
          of the link as a spherical join so that
          \(
            \AltSimplex[\One]\subseteq\TheSphBuild[\One]
          \)
          and
          \(
            \AltSimplex[\Two]\subseteq\TheSphBuild[\Two]
            .
          \)
        \item\label{perp_d}
          The orthogonal projection
          \(
            \OrthProjOf[{\TheSimplex\union\AltSimplex[\One]}]{
              \AltSimplex[\Two]
            }
          \)
          is contained in $\TheSimplex$. (The orthogonal
          projection can be carried out in any Euclidean
          apartment containing $\AltSimplex[\One],\AltSimplex[\Two],$
          and $\TheSimplex$. The result is independent of which
          apartment was chosen.)
      \end{enumerate}
    \end{lemma}
    \begin{proof*}
      \begin{description}
        \item[\pref{perp_b}$\Longrightarrow$\pref{perp_c}]
          This follows from Lemma~\pref{factors}.
        \item[\pref{perp_c}$\Longrightarrow$\pref{perp_d}]
          clear.
        \item[\pref{perp_d}$\Longrightarrow$\pref{perp_b}]
          clear.\qed
      \end{description}
    \end{proof*}

    \begin{lemma}\label{min_face}
      For any horizontal simplex $\TheSimplex$, there is a
      unique face $\TheSimplex[][\min]\faceof\TheSimplex$ such
      that for any proper face $\TheFace\strictfaceof\TheSimplex$,
      we have the equivalence
      \[
        \TheSimplex\setminus\TheFace
        \in
        \HorLinkOf{\TheFace}
        \quad\Longleftrightarrow\quad
        \TheSimplex[][\min]\faceof\TheFace\strictfaceof\TheSimplex
      \]
      More precisely, using any chamber $\TheAffChamber$ containing
      $\TheSimplex$, the face $\TheSimplex[][\min]$ can be described as
      the smallest face of $\TheSimplex$ containing the set
      \[
        \SetOf[{
          \OrthProjOf[\TheSimplex]{\TheAffVertex}
        }]{
          \TheAffVertex\text{\ vertex in\ }\TheAffChamber,\,\,
          \TheHeightOf{\TheAffVertex}\neq\TheHeightOf{\TheSimplex}
        }
        .
      \]
    \end{lemma}
    \begin{proof}
      Note that uniqueness of $\TheSimplex[][\min]$ is obvious.
      It remains to show that for any choice of the chamber
      $\TheAffChamber$ the face $\TheSimplex[][\min]$ defined
      above satisfies
      \[
        \TheSimplex\setminus\TheFace
        \in
        \HorLinkOf{\TheFace}
        \quad\Longleftrightarrow\quad
        \TheSimplex[][\min]\faceof\TheFace\strictfaceof\TheSimplex
      \]
      for each proper face $\TheFace\strictfaceof\TheSimplex$.

      By Lemma~\pref{is_equ}, we have
      \(
        \TheSimplex\setminus\TheFace
        \subseteq
        \HorLinkOf{\TheFace}
      \)
      if and only if every non-equatorial
      vertex $\TheAffVertex\in\TheAffChamber$
      has distance $\PiHalf$ to $\TheSimplex\setminus\TheFace$,
      which by Lemma~\pref{perpendicular} happens if and only if
      \(
        \OrthProjOf[\TheSimplex]{\TheAffVertex}
        \in
        \TheFace
        .
      \)
      This, in turn, is equivalent to
      \(
        \TheSimplex[][\min]\faceof\TheFace
      \)
      by construction of $\TheSimplex[][\min]$.
    \end{proof}
    Let $\TheFace$ be a face of $\TheSimplex$. Note that
    in any Euclidean apartment containing $\TheSimplex$, the
    orthogonal projection onto the affine subspace spanned
    by $\TheFace$ factors through the orthogonal projection
    onto the subspace spanned by $\TheSimplex$. It is now easy
    to make the following:
    \begin{observation}\label{intermediate}
      Suppose $\TheSimplex[][\min]\faceof\TheFace\faceof\TheSimplex$
      for some horizontal simplices.
      Then $\TheSimplex[][\min]=\TheFace[][\min]$.
      In particular,
      \(
        \TheSimplex[][\min]=
        \Parentheses[][\min]{
          \TheSimplex[][\min]
        }
        .
      \)
    \end{observation}
    \begin{proof}
      Let $\TheAffChamber$ be a chamber containing $\TheSimplex$.
      Then, for any vertex $\TheAffVertex\in\TheAffChamber$ not on
      the level of $\TheSimplex$, we have
      \[
        \OrthProjOf[\TheFace]{\TheAffVertex}
        =
        \OrthProjOf[\TheFace]{
          \OrthProjOf[\TheSimplex]{\TheAffVertex}
        }
        =
        \OrthProjOf[\TheSimplex]{\TheAffVertex}
      \]
      since
      \(
        \OrthProjOf[\TheSimplex]{\TheAffVertex}
        \in
        \TheSimplex[][\min]\subseteq\TheFace
      \)
      by hypothesis.
    \end{proof}

    We now define two relations on horizontal simplices.
    We define \notion{going up} as
    \[
      \TheFace \uprel \TheSimplex
      \quad:\Longleftrightarrow\quad
      \TheFace = \TheSimplex[][\min] \neq \TheSimplex
    \]
    and \notion{going down} as
    \[
      \TheSimplex \downrel \TheFace
      \quad:\Longleftrightarrow\quad
      \TheSimplex[][\min] \not\faceof \TheFace \strictfaceof \TheSimplex
      .
    \]
    We define a \notion{move} as either going up or going
    down and write
    \(
      \TheSimplex[\One] \moveto \TheSimplex[\Two]
    \)
    if there is a move from $\TheSimplex[\One]$ to $\TheSimplex[\Two]$.
    The main result of this section is the following
    \begin{prop}\label{uniform_bound}
      There is a uniform bound, depending only on the
      dimension of $\TheAffBuild$, on the length of any
      sequence of moves.
    \end{prop}
    Thus, we can define the \notion{depth}
    \(
      \TheDepthOf{\TheSimplex}
    \)
    of a simplex as the length of a longest sequence of
    moves starting at $\TheSimplex$.
    Assuming for a moment that the depth is well defined,
    we have the following:
    \begin{observation}\label{depth_dec}
      If there is a move from $\TheSimplex[\One]$ to
      $\TheSimplex[\Two]$, then
      \(
        \TheDepthOf{\TheSimplex[\One]} >
        \TheDepthOf{\TheSimplex[\Two]}
      \)
      since we can put the move from
      $\TheSimplex[\One]$ to $\TheSimplex[\Two]$
      in front of a sequence starting at
      $\TheSimplex[\Two]$ and obtain a longer chain
      starting at $\TheSimplex[\One]$.\qed
    \end{observation}
    The remainder of this section is entirely devoted to the
    proof of Proposition~\pref{uniform_bound} and independent
    of the other parts of the paper.

    Let us begin by collecting some elementary properties of
    the two types of moves. We begin with transitivity.
    \begin{lemma}\label{transitivity_up}
      It never happens that
      \(
        \TheSimplex[\One]\uprel\TheSimplex[\Two]\uprel\TheSimplex[\Three]
        .
      \)
      In particular, the symmetric closure of $\uprel$ is
      transitive for silly reasons.
    \end{lemma}
    \begin{proof}
      Suppose
      \(
        \TheSimplex[\One]\uprel\TheSimplex[\Two]\uprel\TheSimplex[\Three]
        .
      \)
      Then, by Observation~\pref{intermediate},
      \(
        \TheSimplex[\One]
        =
        \TheSimplex[\Two][\min]
        =
        \Parentheses[][\min]{
          \TheSimplex[\Three][\min]
        }
        =
        \TheSimplex[\Three][\min]
        =
        \TheSimplex[\Two]
      \)
      contradicting
      \(
        \TheSimplex[\One]\strictfaceof\TheSimplex[\Two]
        .
      \)
    \end{proof}
    \begin{lemma}\label{transitivity_down}
      The relation $\downrel$ is transitive.
    \end{lemma}
    \begin{proof}
      Suppose
      \(
        \TheSimplex[\One]\downrel
        \TheSimplex[\Two]\downrel
        \TheSimplex[\Three]
        .
      \)
      Then
      \(
        \TheSimplex[\One][\min]\not\faceof\TheSimplex[\Two]
        \strictfaceof\TheSimplex[\One]
      \)
      and
      \(
        \TheSimplex[\Two][\min]\not\faceof\TheSimplex[\Three]
        \strictfaceof\TheSimplex[\Two]
        .
      \)
      It follows immediately that
      \(
        \TheSimplex[\Three]\strictfaceof\TheSimplex[\One]
        .
      \)
      Also, $\TheSimplex[\One][\min]\not\faceof\TheSimplex[\Two]$
      and $\TheSimplex[\Three]\faceof\TheSimplex[\Two]$ imply
      that $\TheSimplex[\One][\min]\not\faceof\TheSimplex[\Three]$.
      Thus,
      \(
        \TheSimplex[\One]\downrel\TheSimplex[\Three]
        .
      \)
    \end{proof}

    The next batch of lemmata deals with chains of simplices
    \[
        \TheFace[\One]\uprel\TheSimplex[\One]
        \downrel
        \TheFace[\Two]\uprel\TheSimplex[\Two]
        \downrel\cdots
   \]
    alternatingly going up and down.

    \begin{lemma}\label{simplify}
      If some horizontal simplices satisfy
      \[
        \TheFace[\One]\uprel\TheSimplex[\One]
        \downrel
        \TheFace[\Two]
        ,
      \]
      then we have
      \[
        \TheFace[\One]
        =
        \Parentheses[][\min]{
          \TheFace[\One]\union\TheFace[\Two]
        }
        \text{\ and\ }
        \TheFace[\One]\union\TheFace[\Two]
        \downrel
        \TheFace[\Two]
        .
      \]
      In particular, we have
      \(
        \TheFace[\One]\uprel\TheFace[\One]\union\TheFace[\Two]
        \downrel\TheFace[\Two]
      \)
      unless
      \(
        \TheFace[\One]\downrel\TheFace[\Two]
        .
      \)
    \end{lemma}
    \begin{proof}
      From Observation~\ref{intermediate}, we deduce
      \(
        \TheFace[\One]=\Parentheses[][\min]{
          \TheFace[\One]\union\TheFace[\Two]
        }
        .
      \)
      On the other hand,
      \(
        \TheFace[\One]=\TheSimplex[\One][\min]\not\faceof\TheFace[\Two]
      \),
      whence
      \(
        \TheFace[\Two]\strictfaceof\TheFace[\One]\union\TheFace[\Two]
      \)
      and
      \(
        \TheFace[\One]\union\TheFace[\Two]\downrel\TheFace[\Two].
      \)
    \end{proof}

    \begin{lemma}\label{carrier}
      Let $\TheFace[\One]$ and $\TheFace[\Two]$ be two
      simplices whose union is a horizontal simplex.
      Let $\TheAffVertex\in\LinkOf{\TheFace[\One]\union\TheFace[\Two]}$
      be a vertex,
      and let $\TheFace\subseteq\TheFace[\One]\union\TheFace[\Two]$
      be the \notion{carrier} of
      \(
        \OrthProjOf[{\TheFace[\One]\union\TheFace[\Two]}]{
          \TheAffVertex
        }
        ,
      \)
      i.e., the smallest face of
      \(
        \TheFace[\One]\union\TheFace[\Two]
       \)
      containing
      \(
        \OrthProjOf[{\TheFace[\One]\union\TheFace[\Two]}]{
          \TheAffVertex
        }
        .
      \)
      Then
      \(
        \OrthProjOf[{
          \TheFace[\Two]\union\SetOf{\AltAffVertex}
        }]{
          \TheAffVertex
        }
        \not\in
        \TheFace[\Two]
      \)
      for every vertex $\AltAffVertex\in\TheFace\setminus\TheFace[\Two]$.
      In particular, $\TheAffVertex$ and $\AltAffVertex$
      are in the same irreducible factor of
      $\LinkOf{\TheFace[\Two]}$.
    \end{lemma}
    \begin{proof}
      The point
      \(
        \OrthProjOf[{
          \TheFace[\Two]\union\TheFace[\One]
        }]{
          \TheAffVertex
        }
      \)
      is a convex combination of the vertices in
      \(
        \TheFace[\Two]\union\TheFace[\One]
        .
      \)
      Since $\AltAffVertex$ lies in the carrier of
      \(
        \OrthProjOf[{
          \TheFace[\Two]\union\TheFace[\One]
        }]{
          \TheAffVertex
        },
      \)
      we can infer that the $\AltAffVertex$-coordinate of
      \(
        \OrthProjOf[{
          \TheFace[\Two]\union\TheFace[\One]
        }]{
          \TheAffVertex
        }
      \)
      is non-zero. From
      \(
        \OrthProjOf[{
          \TheFace[\Two]\union\SetOf{\AltAffVertex}
        }]{
          \AltAffVertex
        }
        =
        \AltAffVertex
      \)
      we can now deduce that
      \(
        \OrthProjOf[{
          \TheFace[\Two]\union\SetOf{\AltAffVertex}
        }]{
          \TheAffVertex
        }
        =
        \OrthProjOf[{
          \TheFace[\Two]\union\SetOf{\AltAffVertex}
        }]{
          \OrthProjOf[{
            \TheFace[\Two]\union\TheFace[\One]
          }]{
            \TheAffVertex
          }
        }
      \)
      still has a non-zero $\AltAffVertex$-coordinate, whence
      it cannot lie in $\TheFace[\Two]$.
    \end{proof}
    \begin{lemma}\label{vertical}
      Let $\TheFace[\One]$ and $\TheFace[\Two]$ be two
      simplices whose union is a horizontal simplex.
      Suppose
      \(
        \TheFace[\One]=\Parentheses[][\min]{
          \TheFace[\One]\union\TheFace[\Two]
        }
        .
      \)
      Then any vertex
      \(
        \AltAffVertex\in\TheFace[\One]\setminus\TheFace[\Two]
      \)
      lies in
      $\VertLinkOf{\TheFace[\Two]}$.
    \end{lemma}
    \begin{proof}
      Note that
      \(
        \AltAffVertex\in\TheFace[\One]=\Parentheses[][\min]{
          \TheFace[\One]\union\TheFace[\Two]
        }
        .
      \)
      It follows that there is a vertex
      \(
        \TheAffVertex
        \in
        \LinkOf{\TheFace[\One]\union\TheFace[\Two]}
      \)
      with
      \(
        \TheHeightOf{\TheAffVertex}
        \neq
        \TheHeightOf{\TheFace[\One]\union\TheFace[\Two]}
      \)
      that is a witness for
      \(
        \AltAffVertex\in\TheFace[\One]=\Parentheses[][\min]{
          \TheFace[\One]\union\TheFace[\Two]
        }
        ,
      \)
      i.e., $\AltAffVertex$ belongs to the smallest
      simplex containing
      \(
        \OrthProjOf[{
          \TheFace[\One]\union\TheFace[\Two]
        }]{
          \TheAffVertex
        }
        .
      \)
      It then follows from Lemma~\pref{carrier} that
      \(
        \AltAffVertex
      \)
      belongs to the same irreducible factor of
      $\LinkOf{\TheFace[\Two]}$ as $\TheAffVertex$.
      As
      \(
        \TheAffVertex\in\VertLinkOf{\TheFace[\Two]}
        ,
      \)
      we have
      \(
        \AltAffVertex\in\VertLinkOf{\TheFace[\Two]}
        .
      \)
    \end{proof}
    \begin{cor}
      Let
      \[
        \TheFace[\One]\uprel\TheSimplex[\One]
        \downrel
        \TheFace[\Two]\uprel\TheSimplex[\Two]
      \]
      be horizontal simplices. Then
      \(
        \TheSimplex[\Two]\union\TheFace[\One]
      \)
      is a horizontal simplex.
    \end{cor}
    \begin{proof}
      Note that $\TheFace[\One]\union\TheFace[\Two]$ is a
      horizontal simplex and
      \(
        \Parentheses{\TheFace[\One]\union\TheFace[\Two]}
        \setminus
        \TheFace[\Two]
      \)
      is a simplex in the vertical link of
      $\TheFace[\Two]$ by Lemma~\pref{vertical}.

      On the other hand, $\TheSimplex[\Two]\setminus\TheFace[\Two]$
      is a simplex in the horizontal link of $\TheFace[\Two]$
      as
      $\TheFace[\Two]\uprel\TheSimplex[\Two]$. It follows
      that
      \(
        \Parentheses{\TheFace[\One]\union\TheFace[\Two]}
        \setminus
        \TheFace[\Two]
       \)
      and
      \(
        \TheSimplex[\Two]\setminus\TheFace[\Two]
      \)
      span a simplex in $\LinkOf{\TheFace[\Two]}$. The claim
      follows.
    \end{proof}

    Now, we are ready to discuss shortening of alternating
    chains and to rule out the existence of cycles. We start by
    ruling out cycles of length $\Two$.
    \begin{observation}\label{no_short_cycles}
      There do not exist horizontal simplices
      $\TheFace$ and $\TheSimplex$ with
      \(
        \TheFace\uprel\TheSimplex
      \)
      and
      \(
        \TheSimplex\downrel\TheFace
      \)
      since
      \(
        \TheFace\uprel\TheSimplex
      \)
      implies $\TheFace=\TheSimplex[][\min]$
      whereas
      \(
        \TheSimplex\downrel\TheFace
      \)
      implies
      \(
        \TheSimplex[][\min]
        \not\faceof
        \TheFace
        .
      \)\qed
    \end{observation}

    \begin{lemma}\label{pre_shorten}
      Given an alternating chain
      \[
        \TheFace[\One]\uprel\TheSimplex[\One]
        \downrel
        \TheFace[\Two]\uprel\TheSimplex[\Two]
        ,
      \]
      we have
      \(
        \Parentheses[][\min]{
          \TheFace[\One]\union\TheSimplex[\Two]
        }
        =
        \TheFace[\One]
        .
      \)
    \end{lemma}
    \begin{proof}
      Let $\TheAffVertex\in\LinkOf{\TheFace[\One]\union\TheSimplex[\Two]}$
      with
      \(
        \TheHeightOf{\TheAffVertex}\neq
        \TheHeightOf{\TheFace[\One]\union\TheSimplex[\Two]}
        .
      \)
      Since
      \(
        \TheFace[\Two]\uprel\TheSimplex[\Two]
        ,
      \)
      we have
      \(
        \TheFace[\One]\union\TheFace[\Two]
        \subseteq
        \TheFace[\One]\union\TheSimplex[\Two]
        ;
      \)
      also,
      \(
        \Parentheses[][\min]{
          \TheFace[\One]\union\TheFace[\Two]
        }
        =
        \TheFace[\One]
      \)
      by Lemma~\pref{simplify}.
      Therefore
      \(
        \OrthProjOf[{
          \TheFace[\One]\union\TheFace[\Two]
        }]{
          \TheAffVertex
        }
        \in
        \TheFace[\One]
        .
      \)
      Also note that
      \(
        \OrthProjOf[{
          \TheSimplex[\Two]
        }]{
          \TheAffVertex
        }
        \in
        \TheFace[\Two]
      \)
      as
      \(
        \TheFace[\Two]=\TheSimplex[\Two][\min].
      \)
      Hence,
      \(
        \OrthProjOf[{
          \TheFace[\One]\union\TheSimplex[\Two]
        }]{
          \TheAffVertex
        }
        \in
        \TheFace[\One]\union\TheFace[\Two]
        .
      \)
      We conclude that
      \(
        \OrthProjOf[{
          \TheFace[\One]\union\TheSimplex[\Two]
        }]{
          \TheAffVertex
        }
        =
        \OrthProjOf[{
          \TheFace[\One]\union\TheFace[\Two]
        }]{
          \OrthProjOf[{
            \TheFace[\One]\union\TheSimplex[\Two]
          }]{
            \TheAffVertex
          }
        }
        =
        \OrthProjOf[{
          \TheFace[\One]\union\TheFace[\Two]
        }]{
          \TheAffVertex
        }
        \in
        \TheFace[\One]
        .
      \)
      Thus, by Lemma~\pref{min_face},
      \[
        \Parentheses[][\min]{
          \TheFace[\One]\union\TheSimplex[\Two]
        }
        \faceof
        \TheFace[\One]
        \faceof
        \TheFace[\One]\union\TheFace[\Two]
        \faceof
        \TheFace[\One]\union\TheSimplex[\Two]
        ,
      \]
      whence
      \[
        \Parentheses[][\min]{
          \TheFace[\One]\union\TheSimplex[\Two]
        }
        =
        \Parentheses[][\min]{
          \TheFace[\One]\union\TheFace[\Two]
        }
        =
        \TheFace[\One]
      \]
      by Observation~\pref{intermediate} and Lemma~\pref{simplify}.
    \end{proof}
    \begin{NewTh}<statement>[Satz]{Shortening Lemma}\label{shorten}
      Any alternating chain
      \[
        \TheFace[\One]\uprel\TheSimplex[\One]
        \downrel
        \TheFace[\Two]\uprel\TheSimplex[\Two]
      \]
      can be shortened to
      \[
        \TheFace[\One]\uprel
        \Parentheses{
          \TheFace[\One]\union\TheSimplex[\Two]
        }
        \downrel
        \TheSimplex[\Two]
        \text{\ or\ }
        \TheFace[\One]
        \uprel
        \TheSimplex[\One]
        \downrel
        \TheSimplex[\Two]
        .
      \]
    \end{NewTh}
    In the second case, one actually has a shorter chain
    \(
      \TheFace[\One]\downrel\TheSimplex[\Two]
      .
    \)
    For technical reasons, however, it is more convenient to keep
    the chain going up initially.
    \begin{proof}
      By Lemma~\pref{pre_shorten}, we have
      \(
        \TheFace[\One]=
        \Parentheses[][\min]{
          \TheFace[\One]\union\TheSimplex[\Two]
        }
        .
      \)
      Also, by Observation~\pref{no_short_cycles}, we have
      \(
        \TheFace[\One]
        \neq\TheFace[\Two]=\TheSimplex[\Two][\min]
      \)
      whence $\TheSimplex[\Two]\neq \TheFace[\One]\union\TheSimplex[\Two]$.
      It follows that
      \(
        \TheFace[\One]\union\TheSimplex[\Two]
        \downrel
        \TheSimplex[\Two]
        .
      \)

      If $\TheFace[\One]\neq\TheFace[\One]\union\TheSimplex[\Two]$,
      we find
      \(
        \TheFace[\One]\uprel
        \Parentheses{
          \TheFace[\One]\union\TheSimplex[\Two]
        }
        \downrel
        \TheSimplex[\Two]
        .
      \)

      If  $\TheFace[\One]=\TheFace[\One]\union\TheSimplex[\Two]$,
      we find
      \(
        \TheSimplex[\Two]\strictfaceof\TheFace[\One]
        \strictfaceof\TheSimplex[\One]
      \)
      and therefore
      \(
        \TheSimplex[\One]\downrel\TheSimplex[\Two]
        .
      \)
    \end{proof}
    As a consequence, we can rule out cycles of arbitrary length.
    \begin{cor}\label{no_cycles}
      No sequence of moves enters a cycle.
    \end{cor}
    \begin{proof}
      Since $\uprel$ and $\downrel$ are both transitive
      (Lemmata~\pref{transitivity_up} and~\pref{transitivity_down}),
      any minimum length cycle has to alternate between
      $\uprel$ and $\downrel$.
      By the Shortening Lemma~\pref{shorten}, a minimum length
      cycle can go up at most once. Thus, a minimum
      length cycle is alternating of length two. This, however,
      is ruled out by Observation~\pref{no_short_cycles}.
    \end{proof}

    \begin{lemma}\label{support}
      Let
      \[
        \TheFace[\One]\uprel\TheSimplex[\One]
        \downrel
        \TheFace[\Two]\uprel\TheSimplex[\Two]
        \downrel
        \cdots
        \downrel
        \TheFace[\TheLastIndex-\One]\uprel\TheSimplex[\TheLastIndex-\One]
        \downrel\TheFace[\TheLastIndex]
      \]
      be an alternating chain of horizontal simplices. Then
      \(
        \TheFace[\One]\union\TheFace[\Two]
        \union\cdots\union
        \TheFace[\TheLastIndex]
      \)
      is a simplex.
    \end{lemma}
    \begin{proof}
      First, we use induction to show that
      \(
        \TheFace[\One]\union\TheFace[\TheLastIndex]
      \)
      is a simplex. The case of a length two chain
      \[
        \TheFace[\One]\uprel\TheSimplex[\One]
        \downrel
        \TheFace[\Two]
      \]
      is obvious. For longer chains, we can use the
      transitivity of $\downrel$ and the
      Shortening Lemma~\pref{shorten} to argue that
      \[
        \TheFace[\One]
        \uprel
        \Parentheses{
          \TheFace[\One]\union\TheSimplex[\Two]
        }
        \downrel
        \TheFace[\Three]\uprel\TheSimplex[\Three]
        \downrel\cdots\downrel
        \TheFace[\TheLastIndex-\One]\uprel\TheSimplex[\TheLastIndex-\One]
        \downrel
        \TheFace[\TheLastIndex]
      \]
      or
      \[
        \TheFace[\One]
        \uprel
        \TheSimplex[\One]
        \downrel
        \TheFace[\Three]\uprel\TheSimplex[\Three]
        \downrel\cdots\downrel
        \TheFace[\TheLastIndex-\One]\uprel\TheSimplex[\TheLastIndex-\One]
        \downrel
        \TheFace[\TheLastIndex]
      \]
      is a shorter alternating chain from $\TheFace[\One]$
      to $\TheFace[\TheLastIndex]$, whence
      \(
        \TheFace[\One]\union\TheFace[\TheLastIndex]
      \)
      is a simplex by induction hypothesis.

      Now, we apply this argument to subsequences
      \[
        \TheFace[\TheIndex]\uprel\TheSimplex[\TheIndex]
        \downrel
        \TheFace[\TheIndex+\One]\uprel\TheSimplex[\TheIndex+\One]
        \downrel
        \cdots
        \downrel
        \TheFace[\AltIndex-\One]\uprel\TheSimplex[\AltIndex-\One]
        \downrel\TheFace[\AltIndex]
      \]
      and find that
      \(
        \TheFace[\TheIndex]\union\TheFace[\AltIndex]
      \)
      is a simplex
      for any two indices, $\TheIndex$ and $\AltIndex$.
      Since the Euclidean building $\TheAffBuild$ is a flag
      complex, it follows that
      \(
        \TheFace[\One]\union\TheFace[\Two]
        \union\cdots\union
        \TheFace[\TheLastIndex]
      \)
      is a simplex.
    \end{proof}
    \begin{proof}[of Proposition~\pref{uniform_bound}]
      By Lemma~\pref{support}, for any strictly alternating chain
      there is a simplex that contains its lower terms (i.e., the
      elements to which the move is going down or from where the
      move is going up). This simplex has at most
      \(
        \Two[][\TheDimOf{\TheAffBuild}+\One]-\One
      \)
      faces.
      Since Corollary~\pref{no_cycles} rules
      out any repetitions in a chain, the length of any strictly
      alternating chain is therefore bounded by
      \(
        \Two (\Two[][\TheDimOf{\TheAffBuild}+\One] - \One) + \One
      \)
      which accounts for a possible move down in the beginning
      and a move up at the end.

      Also note that the longest possible sequences of moves
      going down have
      length $\leq\TheDimOf{\TheAffBuild}$, and there are no
      $\uprel$-chains of length $\Two$ or longer by
      Lemma~\pref{transitivity_up}.

      It follows that
      we can take the uniform upper bound to be
      \(
        \TheDimOf{\TheAffBuild}
        \Parentheses{
          \Two (\Two[][\TheDimOf{\TheAffBuild}+\One] - \One) + \One
        }
        .
      \)
    \end{proof}

  \section{Descending Links: the Irreducible Case}%
    \label{sec:links_irreducible}
    We retain hypotheses, notation, and terminology from the previous
    section. In particular, the Euclidean building $\TheAffBuild$
    is still assumed to be irreducible.

    We subdivide $\TheAffBuild$ as follows. Each
    horizontal simplex is barycentrically subdivided.
    Note that any simplex can be written as the join
    of its maximal horizontal faces. Thus, each simplex has an induced
    subdivision. Also note that the subdivision rule is compatible
    with face relations. Thus, we have defined a subdivision
    of $\TheAffBuild$, which we will denote by
    $\TheDivBuild$.
    Note that the vertices of $\TheDivBuild$ are in
    $\One$-$\One$-correspondence with the horizontal
    simplices of $\TheAffBuild$. We denote
    by $\BarycenterOf{\TheSimplex}$ the vertex in $\TheDivBuild$
    corresponding to the horizontal simplex $\TheSimplex$ in
    $\TheAffBuild$. Simplices in $\TheDivBuild$ correspond to
    sets of chains
    \[
      \SetOf{
        \TheSimplex[\One][\One] \strictfaceof
        \TheSimplex[\Two][\One] \strictfaceof \cdots \strictfaceof
        \TheSimplex[{\TheLastIndex[\One]}][\One]
        =:
        \TheSimplex[][\One]
        ,\,\,\ldots\,\,,
        \TheSimplex[\One][\AltLastIndex] \strictfaceof
        \TheSimplex[\Two][\AltLastIndex] \strictfaceof \cdots \strictfaceof
        \TheSimplex[{\TheLastIndex[\AltLastIndex]}][\AltLastIndex]
        =:
        \TheSimplex[][\AltLastIndex]
      }
    \]
    where
   \(
      \TheSimplex[][\One],\TheSimplex[][\Two],\ldots,
      \TheSimplex[][\AltLastIndex]
    \)
    are horizontal faces (of different $\TheHeight$-heights)
    of a common simplex $\TheSimplex$.
    We infer:
    \begin{observation}\label{link_as_join}
      The link of a vertex
      \(
        \BarycenterOf{\TheSimplex}\in\TheDivBuild
      \)
      decomposes as a join
      \[
        \LinkOf{\BarycenterOf{\TheSimplex}}
        =
        \FaceLinkOf{\BarycenterOf{\TheSimplex}}
        \join
        \CofaceLinkOf{\BarycenterOf{\TheSimplex}}
      \]
      where the \notion{face part}
      $\FaceLinkOf{\BarycenterOf{\TheSimplex}}$
      is the barycentric subdivision of the
      boundary $\BoundaryOf{\TheSimplex}$ and
      the \notion{coface part}
      $\CofaceLinkOf{\BarycenterOf{\TheSimplex}}$
      is
      \(
        \LinkOf{\TheSimplex} \subseteq \TheAffBuild
      \)
      with the induced subdivision.\qed
    \end{observation}

    Observe that $\TheHeight$
    and $\TheDepth$ are well-defined on vertices of $\TheDivBuild$.
    Also, each vertex $\TheDivVertex$ (corresponding to the
    horizontal simplex $\TheSimplex$) has a dimension
    $\TheDimOf{\TheDivVertex} := \TheDimOf{\TheSimplex}$.
    We define the Morse function
    \begin{eqnarray*}
      \TheMorseFct
      \mapcolon
      \TheDivBuild
      &
      \longrightarrow
      &
      \TheReals\crossprod\TheReals
      \\
      \TheDivVertex & \mapsto &
      \Pair{
        \TheHeightOf{\TheDivVertex}
      }{
        \Parentheses{\TheDimOf{\TheAffBuild}+\One}
        \TheDepthOf{\TheDivVertex} +
        \TheDimOf{\TheDivVertex}
      }
    \end{eqnarray*}
    In order to
    meaningfully talk about its sublevel sets, we need to endow
    \(
      \TheReals\crossprod\TheReals
    \)
    with an order relation. We do so by lexicographic order, i.e.,
    \(
      \Parentheses{
        \TheReal[\One],\TheReal[\Two]
      }
      \leq
      \Parentheses{
        \AltReal[\One],\AltReal[\Two]
      }
    \)
    if and only if
    \[
      \TheReal[\One]<\AltReal[\One]
    \]
    or
    \[
      \TheReal[\One]=\AltReal[\One]
      \quad\text{and}\quad
      \TheReal[\Two]\leq\AltReal[\Two]
      .
    \]
    In other words, if $\TheHeight$ decides, we follow that
    decision; but if $\TheHeight$ yields a tie, we use
    $\TheDepth$ to break it; and if $\TheDepth$ still does
    not allow us to make a decision, we resort to $\TheDim$.
    \begin{observation}\label{no_horizontal_edges}
      Let $\TheDivVertex[\One]$ and $\TheDivVertex[\Two]$
      be two vertices in $\TheDivBuild$ (corresponding to
      the horizontal simplices $\TheSimplex[\One]$ and
      $\TheSimplex[\Two]$). Suppose $\TheDivVertex[\One]$ and
      $\TheDivVertex[\Two]$ span an edge. Then either
      \(
        \TheHeightOf{\TheDivVertex[\One]}
        \neq
        \TheHeightOf{\TheDivVertex[\Two]}
      \)
      or $\TheSimplex[\One]$ and $\TheSimplex[\Two]$ are nested,
      i.e., one is a face of the other. Note that
      in this case,
      \(
        \TheDimOf{\TheDivVertex[\One]}
        \neq
        \TheDimOf{\TheDivVertex[\Two]}
        .
      \)
      Consequently, there are no $\TheMorseFct$-horizontal edges
      in $\TheDivBuild$.\qed
    \end{observation}
    \begin{observation}\label{sublevel_eq}
      Note that
      \(
        \Parentheses{\TheDimOf{\TheAffBuild}+\One}
        \TheDepthOf{\TheSimplex}
        +
        \TheDimOf{\TheSimplex}
      \)
      is uniformly bounded from above by a constant, say,
      $\TheConst$.
      Then the sublevel complex in $\TheDivBuild$
      spanned by the vertex set
      \(
        \SetOf[
          \TheDivVertex
        ]{
          \TheMorseFctOf{\TheDivVertex}
          \leq
          \Parentheses{\TheLevel,\TheConst}
        }
      \)
      is a subdivision of the sublevel complex in $\TheAffBuild$
      spanned by the vertex set
      \(
        \SetOf[
          \TheAffVertex
        ]{
          \TheHeightOf{\TheAffVertex}
          \leq
          \TheLevel
        }
        .
      \)
      In particular, both sublevel complexes have the
      same connectivity.\qed
    \end{observation}

    Since we put have an order on the range
    \(
      \TheReals\crossprod\TheReals
    \)
    of the Morse function, we can define
    descending links $\DescLinkOf{\TheDivVertex}$ as usual
    as the part of the link arising from those cells that
    contain $\TheDivVertex$ as their unique highest
    vertex.
    \begin{observation}\label{flag}
      The Euclidean building $\TheAffBuild$ is a flag complex
      and so is the subdivision $\TheDivBuild$. It follows that
      the links $\LinkOf{\BarycenterOf{\TheSimplex}}$ are flag
      complexes, too. The descending
      link $\DescLink{\BarycenterOf{\TheSimplex}}$ is therefore
      the subcomplex spanned by all adjacent vertices in
      $\TheDivBuild$ of strictly smaller $\TheMorseFct$-height.

      It follows that the descending links inherits a
      decomposition from the link as a join
      parts
      \[
        \DescLinkOf{\BarycenterOf{\TheSimplex}}
        =
        \FaceLinkOf[][\downarrow]{\BarycenterOf{\TheSimplex}}
        \join
        \CofaceLinkOf[][\downarrow]{\BarycenterOf{\TheSimplex}}
      \]
      where
      \(
        \FaceLinkOf[][\downarrow]{\BarycenterOf{\TheSimplex}}
        :=
        \DescLinkOf{\BarycenterOf{\TheSimplex}}
        \intersect
        \FaceLinkOf{\BarycenterOf{\TheSimplex}}
      \)
      and
      \(
        \CofaceLinkOf[][\downarrow]{\BarycenterOf{\TheSimplex}}
        :=
        \DescLinkOf{\BarycenterOf{\TheSimplex}}
        \intersect
        \CofaceLinkOf{\BarycenterOf{\TheSimplex}}
         .
      \)\qed
    \end{observation}

    \begin{lemma}\label{desc_link_a}
      Let $\TheSimplex$ be a horizontal simplex with
      $\TheSimplex[][\min] \neq \TheSimplex$. Then
      \(
        \DescLinkOf{\BarycenterOf{\TheSimplex}}
      \)
      is contractible.
    \end{lemma}
    \begin{proof}
      Since the descending link decomposes as a join
      \(
        \DescLinkOf{\BarycenterOf{\TheSimplex}}
        =
        \FaceLinkOf[][\downarrow]{\BarycenterOf{\TheSimplex}}
        \join
        \CofaceLinkOf[][\downarrow]{\BarycenterOf{\TheSimplex}}
        ,
      \)
      it suffices to show that the descending face part
      \(
        \FaceLinkOf[][\downarrow]{\BarycenterOf{\TheSimplex}}
      \)
      is contractible.
      Recall that the face part
      $\FaceLinkOf{\BarycenterOf{\TheSimplex}}$
      is just the barycentric subdivision of the sphere
      $\BoundaryOf{\TheSimplex}$.

      Since $\TheSimplex[][\min] \neq \TheSimplex$, we have
      \(
        \TheSimplex[][\min]\uprel\TheSimplex
      \)
      whence
      \(
        \TheDepthOf{\TheSimplex[][\min]}
        >
        \TheDepthOf{\TheSimplex}
        .
      \)
      Consequently,
      \(
        \FaceLinkOf[][\downarrow]{\BarycenterOf{\TheSimplex}}
      \)
      misses the vertex $\BarycenterOf{\TheSimplex[][\min]}$.

      On the other hand, for any proper face
      $\TheFace\strictfaceof\TheSimplex$ with
      \(
        \TheSimplex[][\min]\not\faceof\TheFace
      \)
      we have
      \(
        \TheSimplex\downrel\TheFace
        ,
      \)
      whence
      \(
        \TheDepthOf{\TheFace}< \TheDepthOf{\TheSimplex}
        ,
      \)
      i.e., the descending face part contains
      all vertices $\BarycenterOf{\TheFace}$ for
      \(
        \TheSimplex[][\min]\not\faceof\TheFace\strictfaceof
        \TheSimplex.
      \)

      Note that we cannot say anything about the depth of
      simplices $\TheFace$ with
      \(
        \TheSimplex[][\min]\strictfaceof\TheFace\strictfaceof\TheSimplex
        .
      \)
      Nonetheless, the information we have is enough to deduce that
      \(
        \FaceLinkOf[][\downarrow]{\BarycenterOf{\TheSimplex}}
      \)
      is homotopy equivalent to a once-punctures sphere
      and hence contractible: Let $\TheDisk$ be
      the subcomplex of
      $\FaceLinkOf{\BarycenterOf{\TheSimplex}}$
      spanned by the set
      \(
        \SetOf[
          \BarycenterOf{\TheFace}
        ]{
          \TheSimplex[][\min]\not\faceof\TheFace
          \strictfaceof\TheSimplex
        }
        .
      \)
      The geometric realization of $\TheDisk$ is the
      sphere $\BoundaryOf{\TheSimplex}$
      with the open star of the simplex
      $\TheSimplex[][\min]$ removed.
      Thus, $\TheDisk$ is a closed ball and hence contractible.
      We have seen that $\TheDisk$ is a subcomplex of
      \(
        \FaceLinkOf[][\downarrow]{\BarycenterOf{\TheSimplex}}
        .
      \)
      Projecting
      away from $\BarycenterOf{\TheSimplex[][\min]}$ defines
      a deformation retraction from
      \(
        \FaceLinkOf[][\downarrow]{\BarycenterOf{\TheSimplex}}
      \)
      onto $\TheDisk$. Hence
      \(
        \FaceLinkOf[][\downarrow]{\BarycenterOf{\TheSimplex}}
      \)
      is contractible.
    \end{proof}

    \begin{lemma}\label{desc_link_b}
      Let $\TheSimplex$ be a
      horizontal simplex $\TheSimplex$
      that satisfies $\TheSimplex=\TheSimplex[][\min]$.
      If $\TheAffBuild$ is thick, then
      \(
        \DescLinkOf{\BarycenterOf{\TheSimplex}}
      \)
      is
      \(
        \Parentheses{
          \TheDimOf{\TheAffBuild}-\Two
        }
      \)-connected.
    \end{lemma}
    \begin{proof}
      Again, we use the decomposition
      \(
        \DescLinkOf{\BarycenterOf{\TheSimplex}}
        =
        \FaceLinkOf[][\downarrow]{\BarycenterOf{\TheSimplex}}
        \join
        \CofaceLinkOf[][\downarrow]{\BarycenterOf{\TheSimplex}}
        .
      \)
      First note that for each proper
      face $\TheFace\strictfaceof\TheSimplex$,
      we have
      \(
        \TheSimplex=\TheSimplex[][\min]
        \not\faceof\TheFace\strictfaceof\TheSimplex
        ,
      \)
      i.e.,
      \(
        \TheSimplex\downrel\TheFace
        .
      \)
      Hence,
      \(
        \TheDepthOf{\TheFace}<\TheDepthOf{\TheSimplex}
        .
      \)
      It follows that
      \(
        \FaceLinkOf[][\downarrow]{\BarycenterOf{\TheSimplex}}
        =
        \FaceLinkOf{\BarycenterOf{\TheSimplex}}
        ,
      \)
      which is homeomorphic to the sphere
      \(
        \BoundaryOf{\TheSimplex}
        .
      \)

      We now have to understand the descending coface part
      \(
        \CofaceLinkOf[][\downarrow]{\BarycenterOf{\TheSimplex}}
        .
      \)
      Recall that
      \(
        \CofaceLinkOf{\BarycenterOf{\TheSimplex}}
      \)
      is just a subdivision of
      \(
        \LinkOf{\TheSimplex}
        =
        \HorLinkOf{\TheSimplex}
        \join
        \VertLinkOf{\TheSimplex}
      \)
      where the subdivision of a simplex
      \(
        \TheFace :=
        \TheFace[\hor]\join\TheFace[\ver]\subset
        \HorLinkOf{\TheSimplex}
        \join
        \VertLinkOf{\TheSimplex}
      \)
      is induced by the barycentric subdivisions of all
      its maximal horizontal faces.
      In particular, if $\TheFace[\ver]$ does not contain any
      equatorial vertices then $\TheFace[\hor]$ is a maximal
      horizontal face and the subdivision of $\TheFace$ is
      given as the join of the subdivisions of $\TheFace[\hor]$
      and $\TheFace[\ver]$.

      In general, $\TheFace[\ver]$ might contain
      equatorial vertices. In that case, the subdivision of
      $\TheFace=\TheFace[\hor]\join\TheFace[\ver]$ might not
      naturally split as a join of a vertical and a horizontal
      part. However, we shall see that under the assumption
      $\TheSimplex=\TheSimplex[][\min]$, the descending
      coface part
      \(
        \CofaceLinkOf[][\downarrow]{\BarycenterOf{\TheSimplex}}
      \)
      does decompose as a join of a horizontal and a vertical
      component: each simplex will factor as described in the
      previous paragraph.

      The key observation is that the descending link does not
      contain any equatorial barycenters from the vertical link
      \(
        \VertLinkOf{\TheSimplex}
        ,
      \)
      provided $\TheSimplex[][\min]=\TheSimplex$.
      To see this, consider a horizontal coface
      $\TheCoface$ of $\TheSimplex$ so that
      $\TheCoface\setminus\TheSimplex$ does not contribute to the
      horizontal link. By Lemma~\pref{min_face}, this means
      \(
        \TheCoface[][\min] \not \faceof \TheSimplex
        ,
      \)
      which implies
      \(
        \TheCoface\downrel\TheSimplex
        ,
      \)
      whence
      \(
        \TheDepthOf{\TheCoface} > \TheDepthOf{\TheSimplex}
        .
      \)
      Thus, $\BarycenterOf{\TheCoface}$ is not in the descending
      link of $\BarycenterOf{\TheSimplex}$.

      With this decomposition of
      \(
        \CofaceLinkOf[][\downarrow]{\BarycenterOf{\TheSimplex}}
        ,
      \)
      we are ready to determine its connectivity.
      The horizontal link $\HorLinkOf{\TheSimplex}$
      (barycentrically subdivided) is fully descending.
      For each horizontal coface $\TheCoface$ with
      \(
        \TheCoface[][\min] \faceof \TheSimplex \strictfaceof \TheCoface
        ,
      \)
      we have
      \(
        \TheCoface[][\min] = \TheSimplex[][\min] = \TheSimplex
      \)
      by Observation~\pref{intermediate}. Thus,
      \(
        \TheSimplex\uprel\TheCoface
        ,
      \)
      whence
      \(
        \TheDepthOf{\TheCoface} < \TheDepthOf{\TheSimplex}
        .
      \)
      Thus,
      \(
        \BarycenterOf{\TheCoface}\in
        \CofaceLinkOf[][\downarrow]{\BarycenterOf{\TheSimplex}}
        .
      \)

      Finally, we consider the vertical part of the descending
      link of $\TheSimplex$. Since we have already seen that
      no equatorial simplices of $\VertLinkOf{\TheSimplex}$
      contribute to $\DescLinkOf{\BarycenterOf{\TheSimplex}}$,
      we see that the Busemann function $\TheHeight$ decides which
      vertices contribute. More precisely, let $\AltSimplex$
      be a horizontal simplex in $\VertLinkOf{\TheSimplex}$, then
      \(
        \BarycenterOf{\AltSimplex}\in
        \CofaceLinkOf[][\downarrow]{\BarycenterOf{\TheSimplex}}
      \)
      if and only if $\TheHeightOf{\AltSimplex}<\TheHeightOf{\TheSimplex}$.
      It follows that
      \(
        \CofaceLinkOf[][\downarrow]{\BarycenterOf{\TheSimplex}}
      \)
      is a subdivision of an open hemisphere complex
      \(
        \TheHemisphereCx
      \)
      in the vertical link
      \(
        \VertLinkOf{\TheSimplex}
      \)
      induced by the Busemann gradient. Since we used that
      gradient to separate the vertical and the horizontal factors
      in
      \(
        \LinkOf{\TheSimplex}
        =
        \HorLinkOf{\TheSimplex}\join\VertLinkOf{\TheSimplex},
      \)
      it follows
      from Theorem~\pref{bernd}
      that the open hemisphere complex
      \(
        \VertLinkOf[][\downarrow]{\BarycenterOf{\TheSimplex}}
      \)
      is spherical of dimension $\TheDimOf{\VertLinkOf{\TheSimplex}}$.

      Thus,
      \[
        \DescLinkOf{\BarycenterOf{\TheSimplex}}
        \isom
        \BoundaryOf{\TheSimplex}
        \join
        \HorLinkOf{\TheSimplex}
        \join
        \TheHemisphereCx
      \]
      is spherical of dimension $\TheDimOf{\TheAffBuild}-\One$,
      hence
      \(
        \Parentheses{
          \TheDimOf{\TheAffBuild}-\Two
        }
      \)-connected.
    \end{proof}

  \section{Descending Links: the General Case}%
    \label{sec:links_general}
    Now, we finally drop the irreducibility hypothesis.
    However,
    we add the assumption of thickness. Let
    \[
      \TheAffBuild
      =
      \TheAffBuild[\One]\crossprod\cdots\crossprod
      \TheAffBuild[\TheLastIndex]
    \]
    be a thick Euclidean building written as a product of
    thick irreducible Euclidean
    buildings $\TheAffBuild[\TheIndex]$.
    Suppose we are given a function
    \begin{eqnarray*}
      \TheHeight \mapcolon \TheAffBuild
      &\longrightarrow&
      \TheReals
      \\
      \TupelOf{
        \TheAffPoint[\One],\ldots,\TheAffPoint[\TheLastIndex]
      }
      &\mapsto&
      \Sum[\TheIndex]{\TheCoefficient[\TheIndex]\TheHeightOf[\TheIndex]{\TheAffPoint[\TheIndex]}}
    \end{eqnarray*}
    as a positive $(\TheCoefficient[\TheIndex]>\Zero)$ linear
    combination of Busemann functions
    \(
      \TheHeight[\TheIndex]
      \mapcolon\TheAffBuild[\TheIndex]\rightarrow\TheReals
      .
    \)
    We subdivide all $\TheAffBuild[\TheIndex]$ as in
    Section~\ref{sec:links_irreducible} and put
    \[
      \TheDivBuild :=
      \TheDivBuild[\One]\crossprod\cdots\crossprod\TheDivBuild[\TheLastIndex]
      .
    \]
    Note that $\TheDivBuild$ is a poly-simplicial complex, i.e.,
    each cell is a product of simplices. In particular, we can regard
    $\TheDivBuild$ as a piecewise Euclidean complex.

    Also, we extend $\TheHeight$ to a Morse function
    \begin{eqnarray*}
      \TheMorseFct \mapcolon \TheDivBuild
      &\longrightarrow&
      \TheReals\crossprod\TheReals
      \\
      \TupelOf{
        \BarycenterOf{\TheSimplex[\One]},\ldots,
        \BarycenterOf{\TheSimplex[\TheLastIndex]}
      }
      &\mapsto&
      \TupelOf{
        \Sum[\TheIndex]{\TheCoefficient[\TheIndex]\TheHeightOf[\TheIndex]{\TheSimplex[\TheIndex]}}
        ,
        \Sum[\TheIndex]{
            \Parentheses{\TheDimOf{\TheAffBuild[\TheIndex]}+\One}
            \TheDepthOf[\TheIndex]{\TheSimplex[\TheIndex]}
            +
            \TheDimOf{\TheSimplex[\TheIndex]}
        }
      }
    \end{eqnarray*}
    Note that edges in
    \(
      \TheDivBuild
      =
      \TheDivBuild[\One]\crossprod\cdots\crossprod\TheDivBuild[\TheLastIndex]
    \)
    always arise from an edge in a single factor. Therefore,
    there are no edges in $\TheDivBuild$ horizontal with respect
    to $\TheMorseFct$.

    The descending link
    \(
      \DescLinkOf{\TheDivVertex}
    \)
    of a vertex $\TheDivVertex\in\TheDivBuild$
    is defined as the subcomplex
    of the link $\LinkOf{\TheDivVertex}$ induced by all those
    poly-simplices $\TheCell$ containing $\TheDivVertex$ as the unique
    point in $\TheCell$ where $\TheMorseFct$ is maximal.
    Since addition
    \[
      \Sum \mapcolon
      \Parentheses[][\TheLastIndex]{
        \TheReals\crossprod\TheReals
      }
      \longrightarrow
      \TheReals\crossprod\TheReals
    \]
    is strictly monotonic in each of the $\TheLastIndex$ arguments,
    we deduce the following two observations, the first of
    which strengthens slightly the statement that there are
    no $\TheMorseFct$-horizontal edges.
    \begin{observation}\label{uniq_top_vertex}
      Every cell in $\TheDivBuild$ has a unique $\TheMorseFct$-highest
      vertex.
    \end{observation}
    The second observation nails the structure of descending links.
    They decompose as joins of descending links taken in the
    factors $\TheDivBuild[\TheIndex]$.
    \begin{observation}\label{desc_links_are_joins}
      For each vertex
      \(
        \TheDivVertex
        =
        \TupelOf{
          \TheDivVertex[\One],\TheDivVertex[\Two],\ldots,
          \TheDivVertex[\TheLastIndex]
        }
        \in\TheDivBuild
        ,
      \)
      we have
      \[
        \DescLinkOf{\TheDivVertex}
        =
        \DescLinkOf{\TheDivVertex[\One]}
        \join
        \DescLinkOf{\TheDivVertex[\Two]}
        \join\cdots\join
        \DescLinkOf{\TheDivVertex[\TheLastIndex]}
        .
      \]
    \end{observation}
    \begin{rem}\label{vanishing_coeffs}
      In the more general case where we allow some of the
      coefficient $\TheCoefficient[\TheIndex]$ to vanish,
      the descending links in $\TheDivBuild$ are joins of
      descending links in those $\TheDivBuild[\TheIndex]$
      where $\TheCoefficient[\TheIndex]\neq\Zero$.
    \end{rem}
    \begin{prop}\label{desc_links}
      For any vertex
      \(
        \TheDivVertex
        =\TupelOf{\TheDivVertex[\One],\ldots,\TheDivVertex[\TheLastIndex]}
        \in\TheDivBuild,
      \)
      the descending link $\DescLinkOf{\TheDivVertex}$ is
      $\Parentheses{\TheDimOf{\TheAffBuild}-\Two}$-connected.
    \end{prop}
    \begin{proof}
      We have
      \(
        \DescLinkOf{\TheDivVertex}
        =
        \DescLinkOf{\TheDivVertex[\One]}
        \join
        \DescLinkOf{\TheDivVertex[\Two]}
        \join\cdots\join
        \DescLinkOf{\TheDivVertex[\TheLastIndex]}
        .
      \)
      By Lemmata~\pref{desc_link_a} and \pref{desc_link_b},
      the factor
      \(
        \DescLinkOf{\TheDivVertex[\TheIndex]}
      \)
      is $\Parentheses{\TheDimOf{\TheAffBuild[\TheIndex]}-\Two}$-connected.
      The claim now follows since the join of an
      $\TheConnectivity$-connected space and an
      $\AltConnectivity$-connected space is
      $\Parentheses{\TheConnectivity+\AltConnectivity+\Two}$-connected.
    \end{proof}

    We shall now adapt Morse theory to our situation.
    In order to do so, we pass to a subdivision once more.
    \begin{lemma}\label{final_subdivision}
      Let $\AltPecComplex$ be a piecewise Euclidean complex
      with a map of its vertices into an ordered set so that
      each cell has a unique highest vertex. Then $\AltPecComplex$
      has a simplicial subdivision that (a) does not introduce
      new vertices, (b) does not change the homotopy type of
      sublevel complexes, and (c) does not change the homotopy type
      of descending links.
    \end{lemma}
    \begin{proof}
      We proceed by induction on skeleta. We do not need to subdivide
      the $\One$-skeleton. So assume that the
      \(
        \Parentheses{\TheDimension-\One}
      \)-skeleton is already subdivided. To subdivide
      an $\TheDimension$-cell,
      cone off the subdivision of its boundary from the unique top
      vertex. It is clear that this subdivision rule does not introduce
      new vertices. Since we used the top vertex as the cone point,
      it also does not change (but subdivides) sublevel
      sets and descending links.
    \end{proof}

    Going back to the special situation at hand,
    recall how in the the proof of
    Proposition~\pref{general_pos}, the main point was to find
    a number $\TheEpsilon$ so that for each level $\TheLevel$,
    \[
      \TheHeightOf[][-\One]{
        \LeftOpenInterval{-\infty}{\TheLevel}
      }
      \monorightarrow
      \TheHeightOf[][-\One]{
        \LeftOpenInterval{-\infty}{\TheLevel+\TheEpsilon}
      }
    \]
    induces isomorphisms in homotopy groups
    $\Homotopy[\TheDimension]$
    for $\TheDimension\leq\TheDimOf{\TheAffBuild}-\Two$.
    Back there, we noted that by Observation~\pref{parallel},
    there is an $\TheEpsilon$ so that
    \(
      \TheHeightOf[][-\One]{
        \ClosedInterval{\TheLevel}{\TheLevel+\TheEpsilon}
      }
    \)
    does not contain complete edges. In the presence of
    horizontal edges, this is blatantly false.
    However, it still follows by the same argument
    that we can choose $\TheEpsilon$, independent of
    $\TheLevel$, so that every edge contained in
    \(
      \TheHeightOf[][-\One]{
        \ClosedInterval{\TheLevel}{\TheLevel+\TheEpsilon}
      }
    \)
    must be horizontal. For this $\TheEpsilon$, we have:
    \begin{lemma}\label{eps_inclusion}
      For each level $\TheLevel\in\TheReals$, the inclusion
      \[
        \TheHeightOf[][-\One]{
          \LeftOpenInterval{-\infty}{\TheLevel}
        }
        \monorightarrow
        \TheHeightOf[][-\One]{
          \LeftOpenInterval{-\infty}{\TheLevel+\TheEpsilon}
        }
      \]
      induces isomorphisms in homotopy groups
      $\Homotopy[\TheDimension]$
      for $\TheDimension\leq\TheDimOf{\TheAffBuild}-\Two$.
    \end{lemma}
    \begin{proof}
      First, we replace $\TheAffBuild$ by its subdivision
      $\TheDivBuild$. Note that this does not affect sublevel
      sets. Second, let $\TheDivBuild[\TheLevel]$ be the
      subcomplex of $\TheDivBuild$ spanned by all
      vertices in the sublevel set
      \(
        \TheHeightOf[][-\One]{
          \LeftOpenInterval{-\infty}{\TheLevel}
        }
      \)
      and note that $\TheDivBuild[\TheLevel]$
      is a deformation retract
      of the sublevel set (by pushing in free faces).
      Thus, it suffices to show that the inclusion
      \(
        \TheDivBuild[\TheLevel]
        \monorightarrow
        \TheDivBuild[\TheLevel+\TheEpsilon]
      \)
      induces isomorphisms in homotopy groups
      $\Homotopy[\TheDimension]$
      for $\TheDimension\leq\TheDimOf{\TheAffBuild}-\Two$.

      We shall use regular Bestvina-Brady Morse theory to accomplish
      this remaining task.
      We use Observation~\pref{uniq_top_vertex} and
      Lemma~\pref{final_subdivision} to simplicially subdivide
      $\TheDivBuild$ without changing descending links or sublevel
      sets and without
      introducing new vertices. Since nothing changed, we will
      keep the notation $\TheDivBuild$.

      Now, define a new height function on the sublevel
      complex $\TheDivBuild[\TheLevel+\TheEpsilon]$
      as follows
      \begin{eqnarray*}
        \AltMorseFct \mapcolon \TheDivBuild[\TheLevel+\TheEpsilon]
        & \longrightarrow & \TheIntegers \subset \TheReals
        \\
        \TheDivVertex & \mapsto &
        \begin{cases}
          -\One & \text{\ if\ } \TheDivVertex\in\TheDivBuild[\TheLevel]
          \\
          \Sum[\TheIndex]{
            \Parentheses{\TheDimOf{\TheAffBuild[\TheIndex]}+\One}
            \TheDepthOf[\TheIndex]{\TheSimplex[\TheIndex]}
            +
            \TheDimOf{\TheSimplex[\TheIndex]}
          }
          & \text{\ otherwise}
        \end{cases}
      \end{eqnarray*}
      Note that every
      vertex
      \(
        \TheDivVertex\in
        \TheDivBuild[\TheLevel+\TheEpsilon]
        \setminus
        \TheDivBuild[\TheLevel]
      \)
      has the
      same descending link with respect to $\AltMorseFct$ as it
      has with respect to $\TheMorseFct$. To see this let
      $\AltDivVertex$ be a vertex in the link of $\TheDivVertex$.
      If $\AltDivVertex$ also belongs to
      \(
        \TheDivBuild[\TheLevel+\TheEpsilon]
        \setminus
        \TheDivBuild[\TheLevel]
      \)
      the edge connecting $\AltDivVertex$ and $\TheDivVertex$ is
      $\TheHeight$-horizontal by our choice of $\TheEpsilon$.
      It follows that whether $\AltDivVertex$ is descending is
      determined by the secondary Morse function.
      If $\AltDivVertex\in\TheDivBuild$, it is
      clearly descending with respect to both Morse functions.

      We put
      \(
        \TheDivBuildOf{\TheNumber}
        :=
        \AltMorseFctOf[][-\One]{
          \ClosedInterval{-\One}{\TheNumber}
        }
        .
      \)
      Note that
      \(
        \TheDivBuild[\TheLevel] = \TheDivBuildOf{-\One}
      \)
      and
      \(
        \TheDivBuild[\TheLevel+\TheEpsilon] =
        \TheDivBuildOf{\TheNumber}
      \)
      for large $\TheNumber$.
      By \cite[Lemma~2.5]{Bestvina.Brady:1997},
      passing from $\TheDivBuildOf{\TheNumber}$
      to $\TheDivBuildOf{\TheNumber+\One}$ changes the
      homotopy type exactly by coning off descending links
      of all vertices $\TheDivVertex$ with
      $\AltMorseFctOf{\TheDivVertex}=\TheNumber+\One$.
      By Proposition~\pref{desc_links},
      descending links
      are $\Parentheses{\TheDimOf{\TheAffBuild}-\Two}$-connected and
      the claim follows.
    \end{proof}
    We obtain the following theorem as an easy corollary:
    \begin{theorem}\label{thm:connectivity}
      Let
      \(
        \TheAffBuild
        =
        \TheAffBuild[\One]\crossprod\cdots\crossprod
        \TheAffBuild[\TheLastIndex]
      \)
      be a thick Euclidean building written as a product of
      irreducible Euclidean buildings. Then the spherical
      building at infinity decomposes as a join
      \(
        \BoundaryOf[][\infty]{\TheAffBuild}
        =
        \BoundaryOf[][\infty]{\TheAffBuild[\One]}
        \join\cdots\join
        \BoundaryOf[][\infty]{\TheAffBuild[\TheLastIndex]}
        .
      \)
      Thus, we can think of points in $\BoundaryOf[][\infty]{\TheAffBuild}$
      as convex linear combinations of points in the
      $\BoundaryOf[][\infty]{\TheAffBuild[\TheIndex]}$.
      Let $\TheInfPoint\in\BoundaryOf[][\infty]{\TheAffBuild}$
      be a point at infinity with non-trivial coordinates in each
      $\BoundaryOf[][\infty]{\TheAffBuild[\TheIndex]}$. Equivalently,
      assume that $\TheInfPoint$ is not contained in any subspace
      \(
        \BoundaryOf[][\infty]{\TheAffBuild[\One]}
        \join\cdots\join
        \BoundaryOf[][\infty]{\TheAffBuild[\TheIndex-\One]}
        \join
        \BoundaryOf[][\infty]{\TheAffBuild[\TheIndex+\One]}
        \join\cdots\join
        \BoundaryOf[][\infty]{\TheAffBuild[\TheLastIndex]}
        .
      \)
      Then complements of horoballs centered at $\TheInfPoint$
      are $\Parentheses{\TheDimOf{\TheAffBuild}-\Two}$-connected.
    \end{theorem}
    \begin{proof}
      The proof is the same as in the final steps of
      proving Proposition~\pref{general_pos}:
      By Lemma~\pref{eps_inclusion}, the inclusion
      \[
        \TheHeightOf[][-\One]{
          \LeftOpenInterval{-\infty}{\TheLevel}
        }
        \monorightarrow
        \TheHeightOf[][-\One]{
          \LeftOpenInterval{-\infty}{\TheLevel+\TheEpsilon}
        }
        \monorightarrow
        \TheHeightOf[][-\One]{
          \LeftOpenInterval{-\infty}{\TheLevel+\Two\TheEpsilon}
        }
        \monorightarrow\cdots\monorightarrow
        \TheAffBuild
      \]
      induces isomorphisms in $\Homotopy[\TheDimension]$ for
      $\TheDimension\leq\TheDimOf{\TheAffBuild}-\Two$.
      Since $\TheAffBuild$ is contractible, it follows that
      sublevel sets are
      \(
        \Parentheses{\TheDimOf{\TheAffBuild}-\Two}
      \)-connected.
    \end{proof}
    Using Remark~\pref{vanishing_coeffs}, the same argument shows:
    \begin{theorem}
      Let
      \(
        \TheHeight = \Sum[\TheIndex]{
          \TheCoefficient[\TheIndex]\TheHeight[\TheIndex]
        }
        \mapcolon \TheAffBuild \rightarrow \TheReals
      \)
      be a non-negative linear combination of Busemann functions.
      Then, level- and sublevel-set of $\TheHeight$ in $\TheAffBuild$
      are
      \(
        \Parentheses{
          \Parentheses{
            \Sum[{\TheCoefficient[\TheIndex]\neq\Zero}]{
              \TheDimOf{\TheAffBuild[\TheIndex]}
            }
          }
          -\Two
        }
      \)-connected.
    \end{theorem}
    We note that the thickness hypothesis
    in our connectivity results derives entirely from the
    use of Theorem~\pref{bernd} via Lemma~\pref{desc_link_b}
    and the assumption could be dropped here if it could be
    removed from Schulz' result.

    We remark that the results of this section can be used
    to determine geometric invariants of actions of certain
    $\ThePrimeSet$-arithmetic groups on their associated
    symmetric spaces
    in the number field case \cite{Rehn:2007}. We also note
    that a generalization of our results
    to $\TheReals$-buildings, which
    arise, e.g., as asymptotic cones, would be of interest.

  \section{Finitness Properties of Rank One Groups}
    We are now in a position to prove Theorem~\ref{thm:sharp}.
    Recall that $\TheField$ is a global function field,
    $\TheGroup$ is a noncommutative, absolutely almost simple
    $\TheField$-group of $\TheField$-rank $\One$, and $\ThePrimeSet$
    is a finite set of pairwise inequivalent valuations on $\TheField$.
    We can restate Theorem~\ref{thm:sharp} as follows:
    \begin{theorem}
      The $\ThePrimeSet$-arithmetic group $\TheGroupOf{\OkaRing}$
      is of type \FType{\SumOfLocalRanks-\One} where
      \(
        \SumOfLocalRanks
        :=
        \Sum[\ThePrime\in\ThePrimeSet]{
          \RankOf[{\TheField[\ThePrime]}]{\TheGroup}
        }
      \)
      is the sum of local ranks or equivalently, the
      dimension of the associated Euclidean building
      \(
        \TheAffBuild :=
        \Product[\ThePrime\in\ThePrimeSet]{
          \TheAffBuild[\ThePrime]
        }
        .
      \)
    \end{theorem}
    \begin{proof}
      Let $\TheHoroCollection$ be a collection of closed horoballs as
      in Theorem~\ref{thm:geometry}. By further increasing the
      parameter $\ThePush$ used to define
      $\TheHoroCollection$, we can choose
      the horoballs so that the distance between
      any two exceeds the diameter of
      the polyhedral cells in $\TheAffBuild$. With this choice,
      every cell of $\TheAffBuild$ meets at most one horoball.
      It follows that we can subdivide $\TheAffBuild$ so that
      all horoballs in $\TheHoroballSet$ become subcomplexes of
      the CW-complex $\TheAffBuild$.

      Let
      \(
        \TheSpine :=
        \TheAffBuild \setminus
        \Union[\TheHoroball\in\TheHoroballSet]{
          \InteriorOf{\TheHoroball}
        }
      \)
      denote the complement of the open horoballs. This is a
      CW-subcomplex of $\TheAffBuild$ containing the horospheres
      \(
        \BoundaryOf{\TheHoroball}
        .
      \)
      By Theorem~\ref{thm:geometry}, the group $\TheGroupOf{\OkaRing}$
      acts cocompactly on $\TheSpine$. Cell stabilizers are finite.
      By \cite[Propositions~1.1 and~3.1]{Brown:1987},
      it suffices to show that
      \(
        \TheSpine
      \)
      is $\Parentheses{\SumOfLocalRanks-\Two}$-connected.

      We have seen in Theorem~\ref{thm:connectivity}, that each horosphere
      $\BoundaryOf{\TheHoroball}$ is
      $\Parentheses{\SumOfLocalRanks-\Two}$-connected.
      Let $\TheCollapse$ denote $\TheSpine$ with the horospheres
      in
      \(
        \SetOf[{
          \BoundaryOf{\TheHoroball}
        }]{
          \TheHoroball\in\TheHoroballSet
        }
      \)
      collapsed.
      Collapsing disjoint $\TheDimension$-connected subcomplexes
      independently does not
      affect homotopy groups in dimensions up to $\TheDimension$.
      It follows that
      \(
        \HomotopyOf[\TheDimension]{\TheSpine}
        =
        \HomotopyOf[\TheDimension]{\TheCollapse}
      \)
      for $\TheDimension\leq\SumOfLocalRanks-\Two$.

      On the other hand,
      \(
        \TheCollapse
      \)
      can also be obtained from $\TheAffBuild$ by collapsing
      independently the horoballs $\TheHoroball\in\TheHoroballSet$.
      Since collapsing contractible subcomplexes does not affect
      the homotopy type, we deduce that $\TheCollapse$ is
      contractible. Hence
      $\TheSpine$ is $\Parentheses{\SumOfLocalRanks-\Two}$-connected.
    \end{proof}


\section*{References}
  \begin{references}%
    \newcommand{\nametie}{\,}%
  \begin{article}{Abel91}{Abels:1991}%
    \au{H.\nametie Abels}
    \ti{Finiteness Properties of certain arithmetic groups in the function field case}
    \lo{Israel J.\ Math.\ 76 (1991)}{113 -- 128}
  \end{article}
  \begin{book}{Abra96}{Abramenko:1996}%
    \au{P.\nametie Abramenko}
    \ti{Twin Buildings and Applications to $S$-Arithmetic Groups}
    \lo{Springer LNM~1641 (1996)}{}
  \end{book}
  \begin{article}{Behr69}{Behr:1969}%
    \au{H.\nametie Behr}
    \ti{{\selectlanguage{german}Endliche Erzeugbarkeit arithmetischer
        Gruppen "u{}ber Funktionenk"o{}rpern}}
    \lo{Inventiones mathematicae 7 (1969)}{1 -- 32}
  \end{article}
  \begin{article}{BeBr97}{Bestvina.Brady:1997}%
    \au{M.\nametie Bestvina, N.\nametie Brady}
    \ti{Morse theory and finiteness properties of groups}
    \lo{Inventiones mathematicae 129 (1997)}{445 -- 470}
  \end{article}
  \begin{article}{Bro87}{Brown:1987}%
    \au{K.S.\nametie Brown}
    \ti{Finiteness Properties of Groups}
    \lo{Journal of Pure and Applied Algebra 44 (1987)}{45 -- 75}
  \end{article}
  \begin{article}{BuGo99}{Bux.Gonzalez:1999}
    \au{K-U.\nametie Bux, C.R.\nametie Gonzalez}
    \ti{The Bestvina-Brady Construction Revisited -- Geometric Computation
        of $\Sigma$-Invariants for Right Angled Artin Groups}
    \lo{Journal of the London Mathematical Society (2) 60 (1999)}{793 -- 801}
  \end{article}
  \begin{article}{BW07}{Bux.Wortman:2007}
    \au{K.-U.\nametie Bux, K.\nametie Wortman}
    \ti{Finiteness properties of arithmetic groups over function fields}
    \lo{Inventiones mathematicae 167 (2007)}{355 -- 378}
  \end{article}
  \begin{article}{Hard69}{Harder:1969}
    \au{G.\nametie Harder}
    \ti{{\selectlanguage{german}Minkowskische Reduktionstheorie "uber Funtionenk"orpern}}
    \lo{Inventiones mathematicae 7 (1969)}{33 -- 54}
  \end{article}
  \begin{book}{Marg91}{Margulis:1991}
    \au{G.A.\nametie Margulis}
    \ti{Discrete Subgroups of Semisimple Lie Groups}
    \lo{Ergebnisse der Mathematik 3.\ Folge, vol.\ 17 (Springer, 1991)}{}
  \end{book}
  \begin{book}{PlRa94}{Platonov.Rapinchuk:1994}
    \au{V.\nametie Platonov, A.\nametie Rapinchuk}
    \ti{Algebraic groups and number theory}
    \lo{Pure and Applied Math.\ 139, Academic Press (Boston, 1994)}{}
  \end{book}
  \begin{book}{Ragh72}{Raghunathan:1972}
    \au{M.S.\nametie Raghunathan}
    \ti{Discrete subgroups of Lie Groups}
    \lo{Ergebnisse der Mathematik und ihrer Grenzgebiete 68, Springer New York-Heidelberg (1972)}{}
  \end{book}
  \begin{thesis}{Rehn07}{Rehn:2007}
    \au{W.H.\nametie Rehn}
    \ti{{\selectlanguage{german}Kontrollierter Zusammenhang
    "uber symmetrischen R"aumen}}
    \lo{PhD thesis (Frankfurt am Main, 2007)}{}
  \end{thesis}
  \begin{thesis}{Schu05}{Schulz:2005}
    \au{B.\nametie Schulz}
    \ti{{\selectlanguage{german}Sph"arische Unterkomplexe sph"arischer Geb"aude}}
    \lo{PhD thesis (Frankfurt am Main, 2005)}{}
  \end{thesis}
  \begin{article}{Spri94}{Springer:1994}
    \au{T.A.\nametie Springer}
    \ti{Reduction theory over global fields}
    \lo{Proc.\ Indian Acad.\ Sci.\ Math.\ Sci.\ 104 (1994)}{207 -- 216}
  \end{article}
  \begin{article}{Stuh80}{Stuhler:1980}%
    \au{U.\nametie Stuhler}
    \ti{Homological properties of certain arithmetic groups in
        the function field case}
    \lo{Inventiones mathatematicae 57 (1980)}{263 -- 281}
  \end{article}
  \begin{article}{????}{????}
    \au{???}
    \ti{???}
    \lo{???}{??? -- ???}
  \end{article}
  \end{references}
\end{document}